\documentclass{article}

\usepackage[english]{babel}
\usepackage[T1]{fontenc}
\usepackage[utf8]{inputenc}
\usepackage{amsfonts}
\usepackage{amsmath}
\usepackage{dsfont}
\usepackage{hyperref}
\usepackage[thmmarks,amsmath]{ntheorem}
\usepackage{amssymb}
\usepackage{geometry}
\usepackage{esint}
\usepackage{enumitem}
\usepackage{pgf,tikz,tkz-tab}
\usepackage{pgfplots}
\usepackage{xcolor}
\usepackage{float}

\numberwithin{equation}{section}

{
\newtheorem{Def}{Definition}[section]
}

{
\newtheorem{Lem}[Def]{\textit{Lemma}}
}

{
\newtheorem{Cor}[Def]{\textit{Corollary}}
}

\newtheorem{Prop}[Def]{\textit{Proposition}}

\newtheorem{Theo}[Def]{Theorem}

{
\theorembodyfont{}
\newtheorem{Rem}[Def]{Remark}
}

{
\theoremstyle{break}
\theoremsymbol{$\square$}
\theorembodyfont{}
\newtheorem*{Dem}{Proof}
}

\allowdisplaybreaks

\begin{document}

\title{Path-connectedness of incompressible Euler solutions}

\author{Philippe Anjolras}

\date{}

\maketitle

\begin{abstract}
We study the incompressible Euler equation and prove that the set of weak solutions is path-connected. More precisely, we construct paths of Hölder regularity $C^{1/2}$, valued in $C^0_{t, loc} L^2_x$ endowed with the strong topology. The main result relies on a convex integration construction adapted from the seminal work of De Lellis and Székelyhidi \cite{DeLellisSzekelyhidioriginal}, extending it to a more broader geometric framework, replacing balls with arbitrary convex compact sets. 
\end{abstract}

\tableofcontents

\section{Introduction}

The incompressible Euler equation :  
\begin{equation}
\left\{ \begin{array}{l}
\partial_t u + \nabla_x \cdot (u \otimes u) + \nabla p = 0 \\
\nabla_x \cdot u = 0 
\end{array} \right. \label{equ-Euler-incomp}
\end{equation} 
describes the movement of a perfect homogeneous fluid, with 
\[ (u, p) : (t, x) \in \mathbb{R} \times \mathbb{R}^n \mapsto (u(t, x), p(t, x)) \in \mathbb{R}^n \times \mathbb{R} \]
In this paper, we always assume $n \geq 2$ is fixed. 

Our goal is to demonstrate the following theorem : 

\begin{Theo} Let $(u_0, u_1) \in C^0_t L^2_x$ be two weak solutions of the incompressible Euler equation \eqref{equ-Euler-incomp}. There exists a path $\gamma : [0, 1] \to C^0_t L^2_x$ valued in the set of weak solutions of the incompressible Euler equation, such that $\gamma(0) = u_0$, $\gamma(1) = u_1$, and having local Hölder continuity $\frac{1}{2}$ in the following sense : 
\[ \forall ~ T > 0, ~ \exists ~ C_T > 0, ~ \forall ~ (s, s') \in [0, 1]^2, \quad \Vert \gamma(s) - \gamma(s') \Vert_{L^{\infty}([-T, T], L^2(\mathbb{R}^n))} \leq C_T |s-s'|^{1/2} \] \label{theo-cheminsCdemi}
\end{Theo}

\begin{Rem} ~

\begin{enumerate} 
\item Even though the continuity of the path $\gamma$ is only local in time, for each fixed $s$, $\gamma(s)$ is globally in time in $C^0_t L^2_x$. Note that $L^2_x$ is endowed here with the strong $L^2$ topology. If both $u_0$ and $u_1$ are globally bounded in time, that is $u_0, u_1 \in L^{\infty}_t L^2_x$, the constructed $\gamma(s)$ can also be chosen in the same space, but it is not necessary to assume more than continuity (and so local boundedness) for the argument to hold. 

More generally, since we get pointwise bounds on $\gamma(s)$ depending on $u_0$ and $u_1$, a small improvement of the argument also implies that, whenever both $u_0$ and $u_1$ are in some Lebesgue space, then the whole path can be chosen to lie in the same Lebesgue space. The locality however is specific to a defect of convergence of approximations of unity in $L^{\infty}$. A global control in time could be obtained in $C^{1/2} L^p_t L^q_x$. 

\item The solutions above are supposed to be defined on the whole time interval $\mathbb{R}$. If $u_0$ and $u_1$ are only defined on a given bounded interval $[t_1, t_2]$, we may apply the same reasoning to get a a path of solutions on $[t_1+\varepsilon, t_2-\varepsilon]$ for an arbitrary $\varepsilon > 0$. 

\item The construction is local in space, so it can be applied the same way for solutions on the torus $\mathbb{T}^n$ for instance (or on any flat manifold without boundary). 
\end{enumerate}
\end{Rem}

The previous theorem relies on the theorem \ref{theo-principal}, stated below, that extends the seminal work of De Lellis and Székelyhidi : 

\begin{Theo}[\cite{DeLellisSzekelyhidioriginal}] Let $\overline{e} : \mathbb{R} \times \mathbb{R}^n \to \mathbb{R}^{+}$ be a function in $C^0_{t, x} \cap L^{\infty}_t L^1_x$. Then there exists a weak solution $u$ of the incompressible Euler equation \eqref{equ-Euler-incomp} such that $u \in L^{\infty}_t L^2_x \cap C^0_t L^2_{w, x}$ and $\frac{|u|^2}{2} = \overline{e}$ almost everywhere. \label{theoDLS}
\end{Theo}

The main observation for our extension is that one can rewrite the identity $\frac{|u|^2}{2} = \overline{e}$ as 
\[ j_{B(0, 1)}(u-0)^2 = 2 \overline{e}(t, x) \]
where $B(0, R)$ is the closed ball of center $0$ and radius $R$, $j_{\mathcal{K}}$ is the gauge of the compact convex $\mathcal{K}$. But neither the choice of $\mathcal{K}$ as a ball or of $0$ as reference point here above are necessary, and in the following we will construct solutions satisfying almost everywhere : 
\[ j_{\mathbf{K}}(u-u_0)^2 = \overline{e} \]
for a given function $u_0$ that will be our base point, a function $\mathbf{K}$ valued in compact convex sets containing $0$, and a scaling (or ``energy'') function $\overline{e}$, under some compatibility hypotheses between $u_0$ and $\sqrt{\overline{e}} \mathbf{K}$ (namely, that $u_0$ lies inside the $\Lambda$-convex hull associated to $\sqrt{\overline{e}} \mathbf{K}$ for the incompressible Euler equation). 

To prove theorem \ref{theo-cheminsCdemi}, one chooses a suitable convex set to construct an intermediate solution $u_{1/2}$ between $u_0$ and $u_1$, and then iterates the process to construct a solution for every dyadic number, before completing the path by continuity. 

\paragraph{Comparison with earlier works} The incompressible Euler equation was introduced by Euler in \cite{Eulerequation1755} and have since been studied extensively, see for instance the book \cite{MajdabookEuler}. To study weak solutions in particular, De Lellis and Székelyhidi \cite{DeLellisSzekelyhidioriginal} introduced a convex integration scheme, showing a strong lack of uniqueness in general. Their approach was inspired by geometric insight from Nash \cite{NashC1embeddings}, Gromov \cite{Gromovbook1986} and Tartar \cite{Tartarcompensatedcompactnesspde, Tartarcompensatedcompactnesscons} and continued pioneering works by Scheffer \cite{Schefferfirstconstr} and Shnirelman \cite{Schnirelmansecondconstr}. De Lellis and Székelyhidi also proved that none of the usual energy criteria to be imposed to select admissible solutions allows to recover uniqueness \cite{DeLellisSzekelyhidi2}. The convex integration method has then been used to construct (many) weak solutions for various PDEs arising in fluid mechanics and understand turbulence effects~: most notably, the resolution of the Onsager conjecture \cite{Onsagerhydro, Kolmogorovturbulence} in \cite{IsettOnsager, IsettOnsager2, BDLSVOnsager}, see also \cite{DeLellisSzekelyhidi3, BDLISAnomalousdissip, DeLellisSzekelyhidi4, BuckmasterOnsagerae, BDLSOnsager}, and the non-uniqueness of weak solutions for the Navier-Stokes equation \cite{BVNS}. For a survey about convex integration results in hydrodynamics, we refer to \cite{BVreview}. The convex integration methods were also adapted to the compressible Euler system \cite{CDLKisentropic}, active scalar equations \cite{CFGporous, ShvActivescalar, SzePorous, BSV_SQG}, ideal magnetohydrodynamics \cite{FLSMHD1, FLSMHD2} and other systems of equations \cite{BLNWild, CMBoussinesq, CSzStationary, KYPerona, Ngeostrophic}. 

Our approach follows closely the initial proof of De Lellis and Székelyhidi, and the structure exposed in \cite{VillaniBourbaki}. However, concerning the topology of the set of weak solutions, the original approach proved their density in a weak topology, while we focus here on insights on the strong topology. It is clear that weak solutions of the incompressible Euler equation \eqref{equ-Euler-incomp} cannot be dense in the strong topology of $L^{\infty}_t L^2_x$, since the set of solutions is closed under strong convergence~: we aim to fill a gap here by describing this set in terms of connectedness and regularity of the associated paths. Note that, by changing the base point around which the construction is based and considering general convex compact sets, the pressure is not anymore given by the simple relation 
\begin{equation} p = -\frac{|u|^2}{n} \label{equ-pressure-simple} \end{equation}
and it is just a $L^{\infty}_t L^1_x$ function (defined up to a constant in space function). By changing the ball to a convex compact set, it also becomes difficult to keep precisely track of the energy density of the solution, which was one the main motivation of the initial construction. In particular, it is unclear whether one can constructs continuous paths joining two, say, energy-decreasing weak solutions by a path consisting only of energy-decreasing solutions (or replacing it by any other energy admissibility criterion). 

\paragraph{Outline of the article} We start in section \ref{section-notations-gauge} by setting some notations and exposing the main abstract theorem \ref{theo-principal}, which is a generalization of theorem \ref{theoDLS}. The rest of the section is dedicated to proving elementary results about gauges of convex compact sets endowed with the Hausdorff distance. In section \ref{section-geometricanal}, we recall the geometric analysis of the system, rephrasing it as a differential inclusion and giving the associated Tartar's wave cone, before computing the $\Lambda$-convex hull for generic convex compact sets. In section \ref{section-localized}, we construct localized oscillating solutions that are the building blocks of the convex integration scheme, with a slightly different potential than the one used by De Lellis and Székelyhidi (due to the absence of the relation \eqref{equ-pressure-simple}), and we construct the building blocks of the theory, which are smooth, localized, very oscillating solutions. In this section, we also explain how to compensate the absence of the usual relation $|u+v_{\varepsilon}|^2 = |u|^2 + |v_{\varepsilon}|^2 + 2 u \cdot v_{\varepsilon}$ in the context of a generic convex set, using uniform convexity to obtain a bound by below. Section \ref{section-convexint} contains the proof of theorem \ref{theo-principal}~: at this point, most of the differences with theorem \ref{theoDLS} have already been tackled and the proof is very similar. For the sake of completeness, proofs are given so to keep this paper self-contained. Finally, section \ref{section-chemins} contains the proof of our main theorem \ref{theo-cheminsCdemi}, applying theorem \ref{theo-principal} repeatedly with suitable convex compact sets. 

\section{Convex integration with generic convex sets} \label{section-notations-gauge}

\subsection{Notations and statement of the main result}

\begin{Def} We denote by $\mathcal{S}^n_0$ the set of symmetric trace-free matrices of size $n \times n$. $I_n$ is the identity matrix of size $n$. We denote by $B(u, R)$ the closed Euclidean ball of center $u \in \mathbb{R}^n$ and radius $R \geq 0$. Let us define 
\[ F(u) = u \otimes u - \frac{|u|^2}{n} I_n \]
which maps $\mathbb{R}^n$ into $\mathcal{S}^n_0$. 

We define the (linear) equation of incompressible Euler subsolutions as : 
\begin{equation} \left\{ \begin{array}{l}
\partial_t u + \nabla_x M + \nabla_x q = 0 \\
\nabla_x \cdot u = 0 
\end{array} \right. \label{equ-Euler-lin} \end{equation}
with unknowns $(u, M, q) : \mathbb{R}^{n+1} \to \mathbb{R}^n \times \mathcal{S}^n_0(\mathbb{R}) \times \mathbb{R}$. 
\end{Def}

\begin{Def} We will denote by $\mathfrak{K}$ the set of convex compact sets containing a neighborhood of $0$. We endow this set (of sets) with the Hausdorff distance, denoted $d_{\mathcal{H}}$.

For $\mathcal{K} \in \mathfrak{K}$, recall that its gauge is defined as 
\[ \forall ~ x \in \mathbb{R}^n, \quad j_{\mathcal{K}}(x) = \inf \{ r > 0, x \in r \mathcal{K} \} \]
Recall that this gauge satisfies all the properties of a norm on $\mathbb{R}^n$, except the absolute homogeneity for negative scalars when $\mathcal{K}$ is not symmetric with respect to $0$. 

We will say that $\mathcal{K}$ is $c$-uniformly convex for some $c > 0$ if 
\[ \forall ~ (x, y) \in \mathbb{R}^n, \quad j_{\mathcal{K}}(x+y)^2 + j_{\mathcal{K}}(x-y)^2 \geq 2 j_{\mathcal{K}}(x)^2 + c j_{\mathcal{K}}(y)^2 \]
We will say that $\mathcal{K}$ is uniformly convex if it is $c$-uniformly convex for some $c > 0$. 

For $\mathcal{K} \in \mathfrak{K}$, let us define
\begin{equation} a_{\mathcal{K}} = \inf \left\{ R^2 - |\overline{u}|^2, ~ \mathcal{K} \subset B(\overline{u}, R) \right\} \label{equdef-aK} \end{equation}
\label{defgenconv}
\end{Def}

Now we state the 

\begin{Theo} Let $(\mathbf{K}, \overline{e}) : \mathbb{R} \times \mathbb{R}^n \to \mathfrak{K} \times \mathbb{R}^{+}$ be continuous, such that $\overline{e} \in L^{\infty}_t L^1_x$ and there exists $0 < r_0 < R_0$ and $c_0 > 0$ such that for every $(t, x) \in \mathbb{R}^{n+1}$, 
\begin{itemize}
\item $B(0, r_0) \subset \mathbf{K}(t, x) \subset B(0, R_0)$ ;
\item $\mathbf{K}(t, x)$ is $c_0$-uniformly convex. 
\end{itemize}

Let $z_0 = (v_0, M_0, q_0) \in L^{\infty}_{loc}(\mathbb{R}, L^2(\mathbb{R}^n, \mathbb{R}^n) \times L^1(\mathbb{R}^n, \mathcal{S}^n_0)) \times L^1_{loc}(\mathbb{R}^{n+1}, \mathbb{R})$ be a solution of \eqref{equ-Euler-lin}. Let us also assume that $z_0 \in C^0_{t, x}$, and that 
\[ \begin{aligned} 
\forall ~ (t, x) \in \mathbb{R}^{n+1}, \quad &F(v_0(t, x)) - M_0(t, x) = 0 &\mbox{ if } \overline{e}(t, x) = 0, \quad \\
\mbox{ or } &F(v_0(t, x)) - M_0(t, x) < \frac{\overline{e}(t, x) a_{\mathbf{K}(t, x)}}{n} I_n &\mbox{ if } \overline{e}(t, x) > 0 
\end{aligned} \]
where the inequality is in the sens of symmetric matrices. 

Then for any $\eta > 0$, there exists a weak solution $v$ of the incompressible Euler equation \eqref{equ-Euler-incomp} such that 
\begin{enumerate}
\item $v \in C(\mathbb{R}, L^2_w(\mathbb{R}^n, \mathbb{R}^n))$ 
\item $j_{\mathbf{K}(t, x)}(v(t, x)-v_0(t, x))^2 = \overline{e}(t, x)$ for every $t \in \mathbb{R}$ and almost every $x \in \mathbb{R}^n$
\item $\sup_{t \in \mathbb{R}} \Vert v(t) - v_0(t) \Vert_{H^{-1}(\mathbb{R}^n)} \leq \eta$ 
\item there exists a sequence $(v_k, p_k, f_k)_k = (v_0 + \widetilde{v}_k, q_0 + \widetilde{q}_k - \frac{|v_k|^2}{n}, \nabla_x \cdot (F(v_k) - M_0) + \widetilde{f}_k )$ where $(\widetilde{v}_k, \widetilde{q}_k, \widetilde{f}_k)_k \in C^{\infty}_c \left( \mathbb{R} \times \mathbb{R}^n, \mathbb{R}^n \times \mathbb{R} \times \mathbb{R}^n \right)^{\mathbb{N}}$ and the $(v_k, p_k, f_k)_k$ are solutions of the forced incompressible Euler equation : 
\[ \left\{ \begin{array}{l}
\partial_t v_k + \nabla_x \cdot (v_k \otimes v_k) + \nabla_x p_k = f_k \\
\nabla_x \cdot v_k = 0 
\end{array} \right. \]
Furthermore, $f_k \underset{k \to \infty}{\overset{\mathcal{D}'}{\longrightarrow}} 0$, $(v_k, p_k) \underset{k \to \infty}{\overset{L^2_{t, loc} L^2_x \times \mathcal{D}'}{\longrightarrow}} (v, p)$. 
\end{enumerate} \label{theo-principal}
\end{Theo}

We follow for this theorem the strategy of De Lellis and Szekelyhidi. 

\subsection{Elementary properties of compact sets and their gauges}

\begin{Def} We define 
\begin{gather*}
r_{max}(\mathcal{K}) = \inf \left\{ R > 0, \mathcal{K} \subset B(0, R) \right\} \\
r_{min}(\mathcal{K}) = \sup \left\{ R \geq 0, B(0, R) \subset \mathcal{K} \right\} \end{gather*}
\end{Def}

\begin{Lem} $r_{min}$ and $r_{max}$ are $1$-Lipschitz continuous and positive valued on $\mathfrak{K}$. \label{lem-continuite-rminmax}
\end{Lem}

\begin{Dem}
By definition, $r_{max}(\mathcal{K}) \geq r_{min}(\mathcal{K}) > 0$ for any $\mathcal{K} \in \mathfrak{K}$. 

Let $\mathcal{K}_1, \mathcal{K}_2 \in \mathfrak{K}$, and denote by $d = d_{\mathcal{H}}(\mathcal{K}_1, \mathcal{K}_2)$, $r_i = r_{min}(\mathcal{K}_i)$, $i = 1, 2$. We then have
\[ B(0, r_1) \subset \mathcal{K}_1 \subset \left\{ x \in \mathbb{R}^n, \mbox{dist}(x, \mathcal{K}_2) \leq d \right\} \]
If $d \geq r_1$, then $r_2 > 0 \geq r_1 - d$. Let us assume then $d < r_1$. The function $j_{\mathcal{K}_2}$ is continuous on $\mathbb{R}^n$, so in particular on $S^{n-1}$. 

Let then $e \in S^{n-1}$ be a maximizer of $j_{\mathcal{K}_2}$ on $S^{n-1}$. This means that $j_{\mathcal{K}_2}(e) = \frac{1}{r_2}$. But $B(0, r_2) \subset \mathcal{K}_2$ is convex and $r_2 e \in \partial \mathcal{K}_2$, so $\mathcal{K}_2 \subset \{ x \in \mathbb{R}^n, x \cdot e \leq r_2 \}$. In particular, the orthogonal projection onto $\mathcal{K}_2$ of a point of the form $t e$ for some $t > r_2$ is always $r_2 e$. 

Set
\[ f(x) = \mbox{dist}(x, \mathcal{K}_2) \]
Then for any $t \geq r_2$ one has
\[ f(te) = t - r_2 \]
But by hypothesis $f(r_1 e) \leq d$, so $r_2 \geq r_1 - d$. 

By symmetry, one also has $r_1 \geq r_2 - d$. We deduce that $r_{min}$ is $1$-Lipschitz continuous on $\mathfrak{K}$. 

Finaly, for $r_{max}$, we have by compactness of $\mathcal{K}_1 \in \mathfrak{K}$ that there exists a $x_1 \in \mathcal{K}_1$ such that $|x_1| = r_{max}(\mathcal{K}_1)$, and by definition of the Hausdorff distance there is a $x_2 \in \mathcal{K}_2$ such that $x_1 \in B(x_2, d)$. Then, $r_{max}(\mathcal{K}_2) \geq r_{max}(\mathcal{K}_1) - d$, and we conclude by symmetry. 
\end{Dem}

\begin{Rem} For any $a \in \mathbb{R}^n$, one has
\[ \frac{|a|}{r_{max}(\mathcal{K})} \leq j_{\mathcal{K}}(a) \leq \frac{|a|}{r_{min}(\mathcal{K})} \] \label{rem-equiv-jauge}
\end{Rem}

\begin{Lem} Let $a \in \mathbb{R}^n$. The function 
\[ \mathcal{K} \in \mathfrak{K} \mapsto j_{\mathcal{K}}(a) \in \mathbb{R}^{+} \]
is locally Lipschitz-continuous. More precisely, 
\[ \left| j_{\mathcal{K}}(a) - j_{\mathcal{K'}}(a) \right| \leq \frac{d_{\mathcal{H}}(\mathcal{K}, \mathcal{K}') |a|}{r_{min}(\mathcal{K}) r_{min}(\mathcal{K}')} \] \label{lem-jauge-Lipschitz} 
\end{Lem}

\begin{Dem}
Let $\mathcal{K}, \mathcal{K}' \in \mathfrak{K}$, and set $d = d_{\mathcal{H}}(\mathcal{K}, \mathcal{K}')$. Then 
\[ a \in j_{\mathcal{K}}(a) \mathcal{K} \subset j_{\mathcal{K}}(a) \left\{ x \in \mathbb{R}^n, ~ \mbox{dist}(x, \mathcal{K}') \leq d \right\}
= \left\{ x \in \mathbb{R}^n, \mbox{dist}(x, j_{\mathcal{K}}(a) \mathcal{K}') \leq d j_{\mathcal{K}}(a) \right\} \]
Therefore, 
\[ a \in j_{\mathcal{K}}(a) \mathcal{K}' + \frac{d j_{\mathcal{K}}(a)}{r_{min}(\mathcal{K}')} \mathcal{K}' \]
This implies
\[ j_{\mathcal{K}'}(a) \leq j_{\mathcal{K}}(a) \left( 1 + \frac{d}{r_{min}(\mathcal{K}')} \right) \leq j_{\mathcal{K}}(a) + \frac{d |a|}{r_{min}(\mathcal{K}') r_{min}(\mathcal{K})} \]

We prove likewise that
\[ j_{\mathcal{K}}(a) \leq j_{\mathcal{K}'}(a) \left( 1 + \frac{d}{r_{min}(\mathcal{K})} \right) \leq j_{\mathcal{K}'(a)} + \frac{d |a|}{r_{min}(\mathcal{K}') r_{min}(\mathcal{K})} \]

But $r_{min}$ is positive, continuous and so locally bounded by below, so we deduce that $j_{\mathcal{K}}(a)$ is locally Lipschitz-continuous. 
\end{Dem}

\begin{Cor} The function 
\[ (a, \mathcal{K}) \in \mathbb{R}^n \times \mathfrak{K} \mapsto j_{\mathcal{K}}(a) \in \mathbb{R}^{+} \]
is locally Lipschitz-continuous. \label{lem-jauge-Lipschitz-tot}
\end{Cor}

\begin{Dem}
Let $(a, \mathcal{K}), (a', \mathcal{K}') \in \mathbb{R}^n \times \mathfrak{K}$. 
\begin{align*}
|j_{\mathcal{K}}(a) - j_{\mathcal{K}'}(a')| &\leq |j_{\mathcal{K}}(a) - j_{\mathcal{K}'}(a)| + |j_{\mathcal{K}'}(a) - j_{\mathcal{K}'}(a')|
\end{align*}
by triangular inequality. But by lemma \ref{lem-jauge-Lipschitz} and the fact that $j_{\mathcal{K}}$ satisfies the triangular inequality, we deduce the Lipschitz-continuity. 
\end{Dem}

\begin{Def} We say that a non-empty compact convex $\mathcal{K} \subset \mathbb{R}^n$ is properly convex if, for any $x \in \partial \mathcal{K}$, there exists $y \in \mathbb{R}^n$ and $R \geq 0$ such that $\mathcal{K} \subset B(y, R)$ and $|x-y| = R$. \end{Def}

\begin{Lem} Let $\mathcal{K} \in \mathfrak{K}$ be uniformly convex. Then $\mathcal{K}$ is properly convex. \label{lempropconv} 
\end{Lem}

\begin{Dem}
Let $\mathcal{K}$ be such a set, $c$ be a constant of uniform convexity (see definition \ref{defgenconv}), and $x \in \partial \mathcal{K}$. In particular, $j_{\mathcal{K}}(x) = 1$. By remark \ref{rem-equiv-jauge}, $j_{\mathcal{K}}^2$ also satisfies 
\[ j_{\mathcal{K}}^2(a+b)^2 + j_{\mathcal{K}}^2(a-b)^2 ~ \geq ~ 2 j_{\mathcal{K}}(a)^2 + \frac{c}{r_{max}(\mathcal{K})^2} |b|^2 \]
for any $a, b \in \mathbb{R}^n$, so we recover the usual notion of strict convexity for $j_{\mathcal{K}}^2$. 

Since $y \mapsto j_{\mathcal{K}}^2(y) - \frac{c}{r_{max}(\mathcal{K})^2} |x-y|^2$ is convex and only takes real values, it admits at least one subdifferential at $x$. Let us denote such a subdifferential $v$. Then : 
\[ \forall ~ y \in \mathbb{R}^n, \quad j_{\mathcal{K}}^2(x+y) - \frac{c}{r_{max}(\mathcal{K})^2} |y|^2 \geq j_{\mathcal{K}}^2(x)^2 + v \cdot y \]

Therefore, $\mathcal{K} = \{ y \in \mathbb{R}^n, ~ j_{\mathcal{K}}(y)^2 \leq 1 \}$ is delimited by the graph of a non-degenerate paraboloid in the direction $v$. But since $\mathcal{K}$ is also a compact set, it is bounded. We deduce that, for $R > 0$ sufficiently large and $v' = \frac{v}{|v|}$, $\mathcal{K}$ is contained in $B(x-Rv', R)$. 

Note that $v$ cannot be $0$ : indeed, if it was the case by contradiction, in the direction $x$, for $\varepsilon \in (-1, 1)$, 
\[ j_{\mathcal{K}}((1+\varepsilon) x)^2 = (1 + \varepsilon)^2 j_{\mathcal{K}}(x)^2 \geq j_{\mathcal{K}}^2(x)^2 + \frac{c}{r_{max}(\mathcal{K})^2} \varepsilon^2 |x|^2 \]
But taking $\varepsilon \to 0$ negative, we get 
\[ 2 \varepsilon j_{\mathcal{K}}(x)^2 \geq O(\varepsilon^2) \]
which is a contradiction. This justifies that we can set $v' = \frac{v}{|v|}$ above. 

This proves that $\mathcal{K}$ is properly convex. 
\end{Dem}

Finally, we give an explicit computation of $a_{\mathcal{K}}$ defined by \eqref{equdef-aK} for some specific compacts used in the proof of theorem \ref{theo-cheminsCdemi}, and prove these compacts satisfy the uniform convexity property of definition \ref{defgenconv}. 

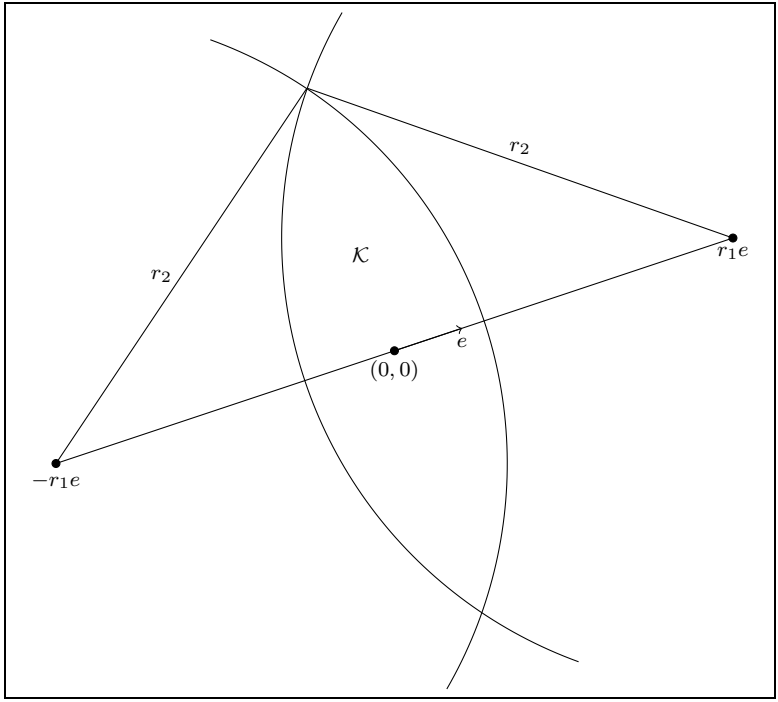
\begin{figure}[H]
\centering
\fbox{
\begin{tikzpicture}[x=3.0cm,y=3.0cm]
\draw[-] (0, 0) -- (1.5, 0.5);
\filldraw[black] (1.5, 0.5) circle (1.5pt) node[below] {\footnotesize $r_1 e$};
\draw[-] (0, 0) -- (-1.5, -0.5); 
\draw[->] (0, 0) -- (0.3, 0.1) node[below] {\footnotesize $e$};
\filldraw[black] (-1.5, -0.5) circle (1.5pt) node[below] {\footnotesize $-r_1e$};
\filldraw[black] (0, 0) circle (1.5pt) node[below] {\footnotesize $(0, 0)$};
\draw plot [variable=\t,samples=100,domain=-30:70] ({-1.5+2.*cos(\t)}, {-0.5+2.*sin(\t)});
\draw plot [variable=\t,samples=100,domain=150:250] ({1.5+2.*cos(\t)}, {0.5+2.*sin(\t)});
\draw (-0.15, 0.5) node[below] {\footnotesize $\mathcal{K}$};
\draw[-] (-1.5, -0.5) -- ({-1.5+2.*cos(56.25)}, {-0.5+2.*sin(56.25)});
\draw ({-1.5+cos(56.25)}, {-0.5+sin(56.25)}) node[left] {\footnotesize $r_2$};
\draw[-] (1.5, 0.5) -- ({1.5+2.*cos(160.6)}, {0.5+2.*sin(160.6)});
\draw ({1.5+cos(160.6)}, {0.5+sin(160.6)}) node[above] {\footnotesize $r_2$};
\end{tikzpicture}
} \caption{\label{figure-cvxbase} The convex of lemma \ref{lem-calcul-aK}}
\end{figure}

\begin{Lem} Let $e \in S^{n-1}$, $r_2 \geq r_1 \geq 0$. Consider the compact (see figure \ref{figure-cvxbase})
\[ \mathcal{K} = B(r_1e, r_2) \cap B(-r_1e, r_2) \]
Then
\[ a_{\mathcal{K}} = r_2^2 - r_1^2 \]
 \label{lem-calcul-aK}
\end{Lem}

\begin{Dem}
$a_{\mathcal{K}}$ being invariant by rotation, we may assume that $e = e_1$ is the first vector of the canonical basis. 

Let $\overline{u}, R$ be such that $\mathcal{K} \subset B(\overline{u}, R)$, and let us write $\overline{u} = (\alpha, \vec \beta) \in \mathbb{R} \times \mathbb{R}^{n-1}$. By symmetry, we assume $\alpha \leq 0$. In this case, the maximal distance between $\overline{u}$ and a point of $\mathcal{K}$ is reached by some point of the form $-r_1 e_1 + r_2 e'$ for some $e' \in S^{n-1}$. Let us write $e' = (\delta, \vec \gamma)$ for a $\delta > 0$ and $\vec \gamma \in e_1^{\perp}$, $|\vec \gamma|^2 + \delta^2 = 1$. In order for $-r_1 e_1 + r_2 e'$ to be in $\mathcal{K}$, it is necessary and sufficient that $e'$ satisfies 
\begin{equation} -r_1 + r_2 e_1 \cdot e' \geq 0 ~ \iff ~ r_2 \delta \geq r_1 \label{equ-preuve-lem-aK} \end{equation}
Let $e'$ such that this inequality holds. 

Then we must have 
\[ |-r_1 e_1 + r_2 e' - \overline{u}| \leq R \]
hence
\[ \begin{aligned}
R^2 - |\overline{u}|^2 &\geq \left| -r_1 e_1 + r_2 e' \right|^2 - 2 \overline{u} \cdot (-r_1 e_1 + r_2 e') \\
&= (r_2 \delta - r_1)^2 + r_2^2 |\vec \gamma|^2 - 2 \alpha (r_2 \delta - r_1) - 2 r_2 \vec \beta \cdot \vec \gamma
\end{aligned} \]
This inequality holds for any $e'$ satisfying \eqref{equ-preuve-lem-aK} : in particular, it is optimal when $\vec \gamma$ is chosen aligned with $\vec \beta$ :  
\begin{align*} 
R^2 - |\overline{u}|^2 &\geq (r_2 \delta - r_1)^2 + r_2^2 (1 - \delta^2) - 2 \alpha (r_2 \delta - r_1) + 2 r_2 |\vec \beta| |\vec \gamma| \\
&\geq (r_2 \delta - r_1)^2 + r_2^2 (1 - \delta^2) - 2 \alpha (r_2 \delta - r_1) \\
&= - 2 r_1 r_2 \delta + r_1^2 + r_2^2 - 2 \alpha (r_2 \delta - r_1) \\
&= -2 \delta r_2 (r_1 + \alpha) + \left( r_1^2 + r_2^2 + 2 \alpha r_1 \right) 
\end{align*}
We then have a degree $1$ polynomial in $\delta$, of main coefficient $-2 r_2 (r_1 + \alpha)$. The optimum in $\delta \in \left[ \frac{r_1}{r_2}, 1 \right]$ (since $e' = (\delta, \vec \gamma)$ is of norm $1$ and $\delta$ satisfies \eqref{equ-preuve-lem-aK}) is then reached for $\delta = 1$ (if $\alpha + r_1 \leq 0$) or $\delta = \frac{r_1}{r_2}$ (if $\alpha + r_1 \geq 0$). In the first case, we then have 
\[ R^2 - |\overline{u}|^2 \geq (r_2 - r_1)^2 - 2 \alpha (r_2 - r_1) \geq (r_2 - r_1)^2 + 2 r_1 (r_2 - r_1) = r_2^2 - r_1^2 \]
And in the second case : 
\[ R^2 - |\overline{u}|^2 \geq \left( r_2 \frac{r_1}{r_2} - r_1 \right)^2 + r_2^2 \left( 1 - \frac{r_1^2}{r_2^2} \right) - 2 \alpha \left( r_2 \frac{r_1}{r_2} - r_1 \right) = r_2^2 - r_1^2 \]
We deduce that $a_{\mathcal{K}} \geq r_2^2 - r_1^2$. 

Finally, considering $\overline{u} = -r_1 e_1$ and $R = r_2$, we have the equality $a_{\mathcal{K}} = r_2^2 - r_1^2$.  
\end{Dem}

\begin{Lem} Let $e \in S^{n-1}$, $a \geq 0$, $R = \sqrt{1 + a^2}$. Set
\[ \mathcal{K} = B(-a e, R) \cap B(ae, R) \in \mathfrak{K} \]
Then $\mathcal{K}$ is $(R-a)^2$-uniformly convex. 

More precisely, if and we set
\[ f(x) = j_{\mathcal{K}}(x)^2 \]
Then, for all $x, y \in \mathbb{R}^n$, it holds that
\begin{equation} f\left( \frac{x+y}{2} \right) \leq \frac{1}{2} \left( f(x) + f(y) \right) - \frac{|x-y|^2}{8} \label{equ-lem-uniforme-convexite-fixe} \end{equation}
 \label{lem-uniforme-convexite-fixe}
\end{Lem}

\begin{Dem}
Let us show that the derivative in the sense of distributions $\nabla^2 f$ satisfies
\[ \nabla^2 f \geq I_n \]
uniformly in $x, a, e$. 

Let us start by computing an explicit formula for $j_{\mathcal{K}}(x)$. Fix $e = e_1$, $x = (x_1, x')$ and, by symmetry, it is enough to consider only the case $x_1 \geq 0$. Then $r = j_{\mathcal{K}}(x)$ is the only non-negative real number such that there exists $e' \in S^{n-1}$ satisfying
\[ -a e_1 + R e' = \frac{x}{r} \]
Therefore, 
\[ R^2 = \left| \frac{x}{r} + a e_1 \right|^2 = \frac{|x|^2}{r^2} + a^2 + 2 a \frac{x_1}{r} ~ \iff ~ (R^2 - a^2) r^2 - 2 a x_1 r - |x|^2 = 0 \]
This is a second order polynomial equation in $r$, whose only positive solution is 
\[ r = \frac{a x_1 + \sqrt{a^2 x_1^2 + (R^2 - a^2) |x|^2}}{R^2 - a^2} \]
Since we assumed $R^2 - a^2 = 1$, this leads to
\[ f(x) = 2 a^2 x_1^2 + |x|^2 + 2 a |x_1| \sqrt{a^2 x_1^2 + |x|^2} \]

Now we derive in the sense of distributions. 
\[ \nabla f(x) = 4 a^2 x_1 e_1 + 2 x + 2 a \epsilon(x_1) e_1 \sqrt{a^2 x_1^2 + |x|^2} + 2 a |x_1| \frac{a^2 x_1 e_1 + x}{\sqrt{a^2 x_1^2 + |x|^2}} \]
where $\epsilon(x_1) = \frac{x_1}{|x_1|}$ is the sign function. Then :  
\[ \begin{aligned}
\nabla^2 f(x) &= 4 a^2 e_1 \otimes e_1 + 2 I_n + 2 a \delta_0(x_1) e_1 \otimes e_1 \sqrt{a^2 x_1^2 + |x|^2} \\
&+ 2 a \epsilon(x_1) \frac{2 a^2 x_1 e_1 \otimes e_1 + x \otimes e_1 + e_1 \otimes x}{\sqrt{a^2 x_1^2 + |x|^2}} \\
&+ 2 a |x_1| \frac{a^2 e_1 \otimes e_1 + I_n}{\sqrt{a^2 x_1^2 + |x|^2}} - 2 a |x_1| \frac{(a^2 x_1 e_1 + x) \otimes (a^2 x_1 e_1 + x)}{(a^2 x_1^2 + |x|^2)^{3/2}} 
\end{aligned} \]
One has $2 a \delta_0(x_1) e_1 \otimes e_1 \geq 0$. 

Then, the third line is also positive semi-definite : 
\[ \begin{aligned}
&2 a |x_1| \frac{a^2 e_1 \otimes e_1 + I_n}{\sqrt{a^2 x_1^2 + |x|^2}} - 2 a |x_1| \frac{(a^2 x_1 e_1 + x) \otimes (a^2 x_1 e_1 + x)}{(a^2 x_1^2 + |x|^2)^{3/2}} \\
&= \frac{2 a |x_1|}{(a^2 x_1^2 + |x|^2)^{3/2}} \left( (a^2 x_1^2 + |x|^2) (a^2 e_1 \otimes e_1 + I_n) - (a^2 x_1 e_1 + x) \otimes (a^2 x_1 e_1 + x) \right) \\
&= \frac{2 a |x_1|}{(a^2 x_1^2 + |x|^2)^{3/2}} \left( |x|^2 I_n - x \otimes x + a^2 \left( x_1^2 I_n + |x|^2 e_1 \otimes e_1 - x_1 e_1 \otimes x - x_1 x \otimes e_1 \right) \right) \\
&\geq \frac{2 a^3 |x_1|}{(a^2 x_1^2 + |x|^2)^{3/2}} \left( x_1^2 \frac{x}{|x|} \otimes \frac{x}{|x|} + |x|^2 e_1 \otimes e_1 - x_1 |x| e_1 \otimes \frac{x}{|x|} - x_1 |x| \frac{x}{|x|} \otimes e_1 \right) \\
&= \frac{2 a^3 |x_1|}{(a^2 x_1^2 + |x|^2)^{3/2}} \left( x_1 \frac{x}{|x|} - |x| e_1 \right) \otimes \left( x_1 \frac{x}{|x|} - |x| e_1 \right) \geq 0 
\end{aligned} \]

Finally, 
\[ \begin{aligned}
&4 a^2 e_1 \otimes e_1 + 2 I_n + 2 a \epsilon(x_1) \frac{2a^2 x_1 e_1 \otimes e_1 + e_1 \otimes x + x \otimes e_1}{\sqrt{a^2 x_1^2 + |x|^2}} \\
&\geq 4 a^2 e_1 \otimes e_1 + 2 I_n + 2 a \epsilon(x_1) \left( e_1 \otimes \frac{x}{\sqrt{a^2 x_1^2 + |x|^2}} + \frac{x}{\sqrt{a^2 x_1^2 + |x|^2}} \otimes e_1 \right) \\
&= 4 a^2 e_1 \otimes e_1 + 2 I_n + \left( 2 a e_1 + \epsilon(x_1) \frac{x}{\sqrt{a^2 x_1^2 + |x|^2}} \right) \otimes \left( 2 a e_1 + \epsilon(x_1) \frac{x}{\sqrt{a^2 x_1^2 + |x|^2}} \right) \\
&- 4 a^2 e_1 \otimes e_1 - \frac{x \otimes x}{a^2 x_1^2 + |x|^2} \\
&\geq I_n 
\end{aligned} \]
as desired. 

This proves \eqref{equ-lem-uniforme-convexite-fixe}, and also 
\[ \forall ~ a, b \in \mathbb{R}^n, \quad f(a+b) + f(a-b) \geq 2 f(a) + |b|^2 \]
On the other hand, by remark \ref{rem-equiv-jauge}, we have that 
\[ \forall ~ b \in \mathbb{R}^n, \quad |b|^2 \geq r_{min}(\mathcal{K})^2 j_{\mathcal{K}}(b)^2 \]
But we can easily compute in this case that 
\[ r_{max}(\mathcal{K}) = \sqrt{R^2 - a^2} = 1, \quad r_{min}(\mathcal{K}) = R - a \]
Therefore, $\mathcal{K}$ is $(R-a)^2$-uniformly convex. 
\end{Dem}

\begin{Rem} Let $\mathcal{K} \in \mathfrak{K}$ be $c$-uniformly convex for some $c > 0$, and let $r > 0$. Then the rescaled convex $r\mathcal{K} \in \mathfrak{K}$ is $c$-uniformly convex. This simply comes from the facte that $j_{r \mathcal{K}} = \frac{1}{r} j_{\mathcal{K}}$. Likewise, $r_{max}(r \mathcal{K}) = r r_{max}(\mathcal{K})$, $r_{min}(r \mathcal{K}) = r r_{min}(\mathcal{K})$. In particular, if we relax the hypothesis of lemma \ref{lem-uniforme-convexite-fixe} to $R > a \geq 0$, without asking $R = \sqrt{1 + a^2}$, we deduce that the associated $\mathcal{K}$ is $\frac{r_{min}(\mathcal{K})^2}{r_{max}(\mathcal{K})^2}$-uniformly convex, with $r_{max}(\mathcal{K}) = \sqrt{R^2 - a^2}$, $r_{min}(\mathcal{K}) = R - a$. \label{rem-unifcvxte} 
\end{Rem}

\begin{Lem} Let $x \in \mathbb{R}^n \mapsto \mathbf{K}(x) \in \mathfrak{K}$ be a continuous function for the Hausdorff distance, such that for every $x$ it holds that
\[ \mathbf{K}(x) = B(-a(x) e(x), R(x)) \cap B(a(x) e(x), R(x)) \]
for some $e(x) \in S^{n-1}$, $a(x) \geq 0$ and $R(x) > \sqrt{2} a(x) \geq 0$. Let us set $r_{max}(x) = r_{max}(\mathbf{K}(x))$. 

Then $\mathbf{K}$ is $\frac{\sqrt{2}-1}{\sqrt{2}+1}$-uniformly convex. 

Let us define the norm
\[ \Vert u \Vert_{L^2_{\mathbf{K}}} := \left( \int_{\mathbb{R}^n} r_{max}(x)^2 j_{\mathbf{K}(x)}(u(x))^2 ~ dx \right)^{1/2} \]
where $u : \mathbb{R}^n \mapsto \mathbb{R}^n$. Then
\begin{equation} \Vert u \Vert_{L^2} \leq \Vert u \Vert_{L^2_{\mathbf{K}}} \leq \sqrt{\frac{\sqrt{2} + 1}{\sqrt{2}-1}} \Vert u \Vert_{L^2} \label{equation-equivnormesconv} \end{equation}

Furthermore, $\Vert \cdot \Vert_{L^2_{\mathbf{K}}}$ is uniformly convex, and this uniformly in $\mathbf{K}$, in the sense that for any $\varepsilon > 0$, there exists a $\delta > 0$ (independent of the choice of $\mathbf{K}$ satisfying the previous hypotheses) such that 
\[ \forall ~ u, v \in L^2(\mathbb{R}^n, \mathbb{R}^n), ~ \mbox{ if } \Vert u \Vert_{L^2_{\mathbf{K}}}, \Vert v \Vert_{L^2_{\mathbf{K}}} \leq 1 \mbox{ and } \left\Vert \frac{u-v}{2} \right\Vert_{L^2_{\mathbf{K}}} > \varepsilon, ~ \mbox{ then } \left\Vert \frac{u+v}{2} \right\Vert_{L^2_{\mathbf{K}}} \leq 1 - \delta \] \label{lem-uniforme-convexite-global}
\end{Lem}

\begin{Dem}
On the one hand it is clear that $|u(x)|^2 \leq r_{max}(x)^2 j_{\mathbf{K}(x)}(u(x))^2$ by definition of $r_{max}$, so we get the left inequality of \eqref{equation-equivnormesconv}. 

On the other hand, 
\[ r_{max}(x)^2 j_{\mathbf{K}(x)}(u(x))^2 \leq \left( \frac{r_{max}(x)}{r_{min}(x)} \right)^2 |u(x)|^2 \]
and so it is enough to show that $\frac{r_{max}}{r_{min}}$ is uniformly bounded. But under the assumption $\sqrt{2} a(x) < R(x)$, we directly have that 
\[ r_{max}(x) = \sqrt{R(x)^2 - a(x)^2}, \quad r_{min}(x) = R(x) - a(x) \]
hence 
\[ 1 \leq \frac{r_{max}(x)}{r_{min}(x)} = \frac{\sqrt{R(x) + a(x)}}{\sqrt{R(x) - a(x)}} \leq \sqrt{\frac{R(x) + \frac{R(x)}{\sqrt{2}}}{R(x) - \frac{R(x)}{\sqrt{2}}}} = \sqrt{\frac{\sqrt{2} + 1}{\sqrt{2} - 1}} \]
which is indeed a constant independent of the choice of $\mathbf{K}$. This ends the proof of \eqref{equation-equivnormesconv}. 

By remark \ref{rem-unifcvxte}, this also proves in this case that $\mathbf{K}$ is $\frac{\sqrt{2}-1}{\sqrt{2}+1}$-uniformly convex. 

Then, for the uniform convexity of $\Vert \cdot \Vert_{L^2_{\mathbf{K}}}$, we notice that we can renormalize 
\[ \mathbf{K}'(x) = \frac{\mathbf{K}(x)}{r_{max}(x)} \]
so that we have $r_{max}(x) = 1$, and $R(x) = \sqrt{1 + a(x)^2}$. This doesn't change the norm : 
\[ \Vert u \Vert_{L^2_{\mathbf{K}}}^2 = \Vert u \Vert_{L^2_{\mathbf{K}'}}^2 = \int j_{\mathbf{K}'(x)}(u(x))^2 ~ dx \]
In what follows, we may therefore assume that $\mathbf{K}(x)$ satisfies $r_{max}(x) = 1$ for every $x$. 

We apply lemma \ref{lem-uniforme-convexite-fixe} and in particular formula \ref{equ-lem-uniforme-convexite-fixe} : if $u, v \in L^2$ satisfy $\Vert u \Vert_{L^2_{\mathbf{K}}}, \Vert v \Vert_{L^2_{\mathbf{K}}} \leq 1$, then 
\[ \begin{aligned}
\left\Vert \frac{u+v}{2} \right\Vert_{L^2_{\mathbf{K}}}^2 
&= \int_{\mathbb{R}^n} j_{\mathbf{K}(x)}\left( \frac{u(x) + v(x)}{2} \right)^2 ~ dx \\
&\leq \int_{\mathbb{R}^n} \frac{1}{2} \left( j_{\mathbf{K}(x)}(u(x))^2 + j_{\mathbf{K}(x)}(v(x))^2 \right) - \frac{|u(x)-v(x)|^2}{8} ~ dx \\
&\leq 1 - \frac{1}{2} \left\Vert \frac{u-v}{2} \right\Vert_{L^2}^2 
\end{aligned} \]
But by the equivalence of norms \eqref{equation-equivnormesconv} already demonstrated, we have for some $C > 0$ (independent of $\mathbf{K}$) that
\[ \left\Vert \frac{u+v}{2} \right\Vert_{L^2_{\mathbf{K}}} \leq \sqrt{1 - C \left\Vert \frac{u-v}{2} \right\Vert_{L^2_{\mathbf{K}}}^2} \]
which show uniform convexity of the norm, and this uniform convexity is uniform in the choice of $\mathbf{K}$. 
\end{Dem}

\section{Geometric analysis} \label{section-geometricanal}

\subsection{Differential inclusion}

\begin{Prop} The following systems of equations are equivalent : on the one hand, the incompressible Euler equation 
\[ \left\{ \begin{array}{l}
\partial_t v + \nabla_x \cdot (v \otimes v) + \nabla_x p = 0 \\
\nabla_x \cdot v = 0 
\end{array} \right. \]
where $(v, p) \in L^{\infty}(\mathbb{R} \times \mathbb{R}^n, \mathbb{R}^n \times \mathbb{R})$, and on the other hand, the equation of subsolutions of incompressible Euler with the inclusion : 
\begin{gather} \left\{ \begin{array}{l}
\partial_t v + \nabla_x \cdot M + \nabla_x q = 0 \\
\nabla_x \cdot v = 0 
\end{array} \right. \\
M = v \otimes v - \frac{|v|^2}{n} I_n ~ \mbox{ almost everywhere} \end{gather}
where $(v, M, q) \in L^{\infty}(\mathbb{R} \times \mathbb{R}^n, \mathbb{R}^n \times \mathcal{S}^n_0 \times \mathbb{R})$. \label{prop-equivalence-systemes}
\end{Prop}

\begin{Dem}
It follows from setting $M = v \otimes v - \frac{|v|^2}{n}$, $q = p + \frac{|v|^2}{n}$. 
\end{Dem}

We then define
\begin{equation} K = \left\{ (v, M, q) \in \mathbb{R}^n \times \mathcal{S}^n_0 \times \mathbb{R}, ~ M = v \otimes v - \frac{|v|^2}{n} I_n \right\} \label{equdef-Kinclusion} \end{equation}

Given a convex compact $\mathcal{K} \in \mathfrak{K}$, we extend the previous definition as
\begin{equation} \overline{\mathcal{K}} = \left\{ (v, M, q) \in K, ~ v \in \partial \mathcal{K} \right\} \label{equdef-Kcvxincl} 
\end{equation}
If also given a base point $z_0 = (v_0, M_0, q_0) \in \mathbb{R}^n \times \mathcal{S}^n_0(\mathbb{R}) \times \mathbb{R}$ and a scaling parameter $r > 0$, we set
\begin{equation} \overline{\mathcal{K}(z_0, r)} = \left\{ (v, M, q) \in K, ~ v - v_0 \in r \partial \mathcal{K} \right\} \label{equdef-Kcvxincl-basescal} 
\end{equation}

\subsection{Tartar's wave cone}

Given a first order linear system $\sum_i A^i \partial_i z = 0$, recall that the set of vectors $a$ such that there exists $\xi \neq 0$ satisfying $\sum_i \xi_i A^i a = 0$, which is equivalent to $x \mapsto a \varphi(\xi \cdot x)$ being a solution of the system for any (say, smooth) function $\varphi$, is called Tartar's wave cone. (The $A^i$ are matrices.) 

\begin{Prop} The wave cone associated to the system \eqref{equ-Euler-lin} of linear subsolutions is
\[ \Lambda = \left\{ (v, M, q) \in \mathbb{R}^n \times \mathcal{S}^n_0 \times \mathbb{R}, ~ \mbox{det} \begin{pmatrix} v & M + q I_n \\ 0 & v^T \end{pmatrix} = 0 \right\} \]
In particular, for any $(v, M) \in \mathbb{R}^n \times \mathcal{S}^n_0$, there exists $q \in \mathbb{R}$ such that $(v, M, q) \in \Lambda$. Furthermore, for any $v_0 \in \mathbb{R}^n$, there exists $p_0 \in \mathbb{R}$ and $\xi \in \mathbb{R}^n \setminus \{ 0 \}$ such that for any function $h \in L^1_{loc}(\mathbb{R})$ and any $\varepsilon > 0$, the couple 
\[ (v, p)(t, x) = (v_0, p_0) h\left( \frac{\xi \cdot x}{\varepsilon} \right) \]
is a solution to the incompressible Euler equation. \label{prop-cone-tartar}
\end{Prop}

\begin{Dem}
Let $U = \begin{pmatrix} v & M + q I_n \\ 0 & v^T \end{pmatrix}$. Then the system \eqref{equ-Euler-lin} can be rewritten as $\nabla_{t,x} \cdot U = 0$. In particular, belonging to the wave cone is equivalent to the existence of a couple $(c, \xi) \in \mathbb{R} \times \mathbb{R}^n$ different from $(0, 0)$ such that $U \begin{pmatrix} c \\ \xi \end{pmatrix} = 0$, which is equivalent to $\mbox{det}(U) = 0$. 

Then, let $(v, M) \in \mathbb{R}^n \times \mathcal{S}^n_0$. Let us first assume $v \neq 0$. Since $n \geq 2$, there exists $\xi \neq 0$ such that $v \cdot \xi = 0$. Therefore, 
\[ U \begin{pmatrix} c \\ \xi \end{pmatrix} = 0 ~ \iff ~ c v + M \xi + q \xi = 0 \]
We then choose $c$ to cancel in the above equation the component along $v$. Then, if we denote by $P_v$ the projection operator on $v^{\perp}$, it is enough to get $P_v M \xi + q \xi = 0$. But $P_v M$ is a symmetric linear operator acting on $v^{\perp}$, which is of dimension at least $1$. So we choose more precisely $\xi$ as an eigenvector of $P_v M$ in $v^{\perp}$, and $q$ the opposite of the associated eigenvalue. 

If $v = 0$, 
\[ U \begin{pmatrix} c \\ \xi \end{pmatrix} = 0 ~ \iff ~ M \xi + q \xi = 0 \]
and we may simply choose $c = 1$ and $\xi = 0$. 

Finally, let $v_0 \in \mathbb{R}^n$. Let us set $M = v_0 \otimes v_0 - \frac{|v_0|^2}{n} I_n$. We notice then that if $\xi \in v^{\perp}$, we also have $M \xi \in v^{\perp}$. In particular, in the previous construction, $c = 0$. Finally, we conclude using the equivalence property \ref{prop-equivalence-systemes}.  
\end{Dem}

\begin{Lem} Let $(a, b) \in \mathbb{R}^n \times \mathbb{R}^n$ be such that $a \neq b$. Set $q = \frac{|b|^2 - |a|^2}{n}$, and 
\[ \lambda = (b-a, F(b) - F(a), q) \]
Then $\lambda \in \Lambda$ and $\lambda \neq 0$. \label{lem-lambda-difference}
\end{Lem}

\begin{Dem}
In dimension at least $2$, there exists a non zero $\xi$ such that $\xi \cdot (a-b) = 0$. Let us set $c = -\xi \cdot a = -\xi \cdot b$. Then
\[ \begin{pmatrix} b-a & F(b) - F(a) + q I_n \\ 0 & (b-a)^T \end{pmatrix} \begin{pmatrix} c \\ \xi \end{pmatrix} = \begin{pmatrix} cb - ca + b (b \cdot \xi) - a (a \cdot \xi) - \frac{|b|^2 - |a|^2}{n} \xi + q \xi \\ (b-a) \cdot \xi \end{pmatrix} = 0 \]
This justifies $\lambda \in \Lambda$ by proposition \ref{prop-cone-tartar}. 

$\lambda \neq 0$ as soon as $a \neq b$. 
\end{Dem}

\begin{Def} Let $\mathcal{K} \in \mathfrak{K}$. We define
\[ \Lambda_{\mathcal{K}} = \left\{ t(b-a, F(b) - F(a), q), ~ t \in \mathbb{R}^{+}, (b, a) \in (\partial \mathcal{K})^2, q = \frac{|b|^2 - |a|^2}{n}, |b|^2 |a|^2 \neq (a \cdot b)^2 \right\} \] 
If given a base point $z_0$ and a scaling parameter $r > 0$, we also define
\[ \Lambda_{\mathcal{K}}(z_0, r) = \left\{ t (b-a, F(b) - F(a), q), ~ t \in \mathbb{R}^{+}, (b, a) \in (v_0 + r\partial \mathcal{K})^2, q = \frac{|b|^2 - |a|^2}{n}, |b|^2 |a|^2 \neq (a \cdot b)^2 \right\} \]
\label{def-Lambda-de-K}
\end{Def}

The condition $|b|^2 |a|^2 \neq (a \cdot b)^2$ simply means that $(a, b)$ are not aligned. Notice that, if $e \in S^{n-1}$, the line $e \mathbb{R}$ intersects $\partial \mathcal{K}$ in exactly two points. The non-alignment condition ensures $a \neq b$, and it will be useful when defining a potential. 

\subsection{Lambda-convex hull of K}

Given a cone $C$, recall that the $C$-convex hull $B$ of a set $A$ is the smallest set $B$ containing $A$ such that, for every couple $(a, b) \in B$, if $b-a \in C$, then $[a, b] \subset B$. 

\begin{Prop} For every uniformly convex $\mathcal{K} \in \mathfrak{K}$, the $\Lambda$-convex hull of $\overline{\mathcal{K}}$ is the convex hull of $\overline{\mathcal{K}}$, and is equal to
\begin{align*} \overline{\mathcal{K}}^{co} = \Biggl\{ (v, M, q) \in \mathcal{K} \times &\mathcal{S}^n_0 \times \mathbb{R}, ~ \forall (\overline{u}, R) \in \mathbb{R}^n \times \mathbb{R}^{+*} \mbox{ such that } \mathcal{K} \subset B(\overline{u}, R), \\
&v \otimes v - M \leq \frac{|v|^2 + R^2 - |v - \overline{u}|^2}{n} I_n \Biggl\} \end{align*}
Furthermore, $\overline{\mathcal{K}} = \overline{\mathcal{K}}^{co} \cap \{ (v, M, q), v \in \partial \mathcal{K} \}$. 

Moreover, there exists a constant $C = C(n) > 0$ such that for every $z = (v, M, q)$ in the interior of $\overline{\mathcal{K}}^{co}$, there exists $\lambda = (\overline{v}, \overline{M}, \overline{q}) \in \Lambda_{\mathcal{K}}$ such that $[z-\lambda, z+\lambda] \in \mbox{int} ~ \overline{\mathcal{K}}^{co}$ and 
\begin{equation} \mbox{dist}\left( \left[ z - \lambda, z + \lambda \right], \partial (\overline{\mathcal{K}}^{co}) \right) \geq \frac{1}{2} d(z, \partial (\overline{\mathcal{K}}^{co})), \quad |\overline{v}| \geq C \mbox{dist}(v, \partial \mathcal{K}) \label{equ-contdist-propenvconv} \end{equation}
In particular, $C$ does not depend on $\mathcal{K}$. \label{prop-enveloppe-convexe}
\end{Prop}

\begin{Dem}
Let us set the notation
\begin{align*} 
C_{\mathcal{K}} := \Biggl\{ (v, M, q) \in \mathcal{K} \times &\mathcal{S}^n_0 \times \mathbb{R}, ~ \forall (\overline{u}, R) \in \mathbb{R}^n \times \mathbb{R}^{+*} \mbox{ such that } \mathcal{K} \subset B(\overline{u}, R), \\
&v \otimes v - M \leq \frac{|v|^2 + R^2 - |v - \overline{u}|^2}{n} I_n \Biggl\} \end{align*}
Note that $\overline{\mathcal{K}} \subset C_{\mathcal{K}}$, with $\overline{\mathcal{K}}$ defined by \eqref{equdef-Kcvxincl}. 

\paragraph{Step 1 : $C_{\mathcal{K}}$ is convex.} Indeed, if we note
\[ \mathcal{B} = \{ (\overline{u}, R) \in \mathbb{R}^n \times \mathbb{R}^{+*}, \mathcal{K} \subset B(\overline{u}, R) \} \]
we have that
\[ \begin{aligned}
C_{\mathcal{K}} &= \left( \mathcal{K} \times \mathcal{S}^n_0 \times \mathbb{R} \right) ~ \cap \bigcap_{(\overline{u}, R) \in \mathcal{B}} \left\{ (v, M, q) \in \mathbb{R}^n \times \mathcal{S}^n_0 \times \mathbb{R}, v \otimes v - M \leq \frac{|v|^2 + R^2 - |v - \overline{u}|^2}{n} I_n \right\} \\
&=: \left( \mathcal{K} \times \mathcal{S}^n_0 \times \mathbb{R} \right) ~ \cap \bigcap_{(\overline{u}, R) \in \mathcal{B}} C_{\mathcal{K}}(\overline{u}, R) 
\end{aligned} \]
The intersection of convex sets is convex, so it is enough to prove that each $C_{\mathcal{K}}(\overline{u}, R)$ is convex. We rewrite the inequality defining $C_{\mathcal{K}}(\overline{u}, R)$ as :  
\[ v \otimes v - M + \frac{2 v \cdot \overline{u}}{n} I_n \leq \frac{R^2 - |\overline{u}|^2}{n} I_n \]
Define then, for $\sigma_{max} : \mathcal{S}_n(\mathbb{R}) \to \mathbb{R}$ the highest eigenvalue function, 
\[ \Phi : (v, M) \mapsto \sigma_{max}\left( v \otimes v - M + \frac{2 v \cdot \overline{u}}{n} I_n \right) := \sup_{\xi \in \mathbb{R}^n, |\xi| = 1} \left( (v \cdot \xi)^2 - \xi \cdot M\xi + \frac{2 v \cdot \overline{u}}{n} \right) \]
$\Phi$ is the supremum of convex functions, so it is convex. Therefore, 
\[ C_{\mathcal{K}}(\overline{u}, R) = \Phi^{-1}\left( \frac{R^2 - |\overline{u}|^2}{n} \right) \]
is convex. So $C_{\mathcal{K}}$ is convex. 

\paragraph{Step 2 : $\overline{\mathcal{K}}$ contains all extremal points of $C_{\mathcal{K}}$.} Let $(v, M, q) \in C_{\mathcal{K}}$. Denote by $\lambda_1 \leq ... \leq \lambda_n$ the list of ordered eigenvalues of $v \otimes v - M$. Let also $(e_i)$ be an orthonormal basis of associated eigenvectors. 

\begin{itemize} \item If $v \in \partial \mathcal{K}$, since $\mathcal{K}$ is properly convex by lemma \ref{lempropconv}, there exists $(\overline{u}, R) \in \mathcal{B}$ such that $v \in \partial B(\overline{u}, R)$. But since $(v, M, q) \in C_{\mathcal{K}}$, this implies 
\[ 0 \leq M - v \otimes v + \frac{|v|^2}{n} I_n \]
The matrix on the right is trace-free, so it has to be the zero matrix. In this case, we deduce that $(v, M) \in \overline{\mathcal{K}}$. 

\item If there exists $(\overline{u}, R) \in \mathcal{B}$ such that $\lambda_1 = \frac{|v|^2 + R^2 - |v - \overline{u}|^2}{n}$ then since $(v, M, q) \in C_{\mathcal{K}}$ we deduce that all the eigenvalues have to be equal to $\lambda_1$ and 
\[ v \otimes v - M = \frac{|v|^2 + R^2 - |v - \overline{u}|^2}{n} I_n \]
By taking the trace, we find that 
\[ R^2 = |v - \overline{u}|^2 \]
so $v \in \partial B(\overline{u}, R)$, hence $v \in \partial \mathcal{K}$ and we recover the previous case. 

\item Assume now that $v \notin \partial \mathcal{K}$ and that for every $(\overline{u}, R) \in \mathcal{B}$ we have $\lambda_1 < \frac{|v|^2 + R^2 - |v - \overline{u}|^2}{n}$ (and in particular $v \notin \partial B(\overline{u}, R)$). 

We start by proving a quantitative estimate on the upper bound on $\lambda_1$. For this, let us show that 
\[ \inf \left\{ R^2 - |v - \overline{u}|^2, ~ (\overline{u}, R) \in \mathcal{B} \right\} \]
is reached at a point $(\overline{u}, R)$ (recall that all the balls are closed). For this, consider a minimizing sequence $(\overline{u}_m, R_m)$. 

If $(R_m)$ is bounded, then $(\overline{u}_m)$ also (because $\mathcal{K}$ is compact so bounded), so up to extracting a subsequence we may assume that both $(R_m)$ and $(\overline{u}_m)$ converge to $R_{\infty}$ and $\overline{u}_{\infty}$ respectively. Then, for every $u \in \mathcal{K}$, $u \in B(\overline{u}_m, R_m)$ for every $m \in \mathbb{N}$, which means $|u - \overline{u}_m| \leq R_m$. The function $f_u : (\overline{u}, R) \mapsto R - |u - \overline{u}|$ is continuous and non-negative evaluated in the $(\overline{u}_m, R_m)$, so by passing to the limit it is also non-negative in $(\overline{u}_{\infty}, R_{\infty})$, so $u \in B(\overline{u}_{\infty}, R_{\infty})$ and therefore $\mathcal{K} \subset B(\overline{u}_{\infty}, R_{\infty})$. On the other hand, it is clear that 
\[ R_{\infty}^2 - |v - \overline{u}_{\infty}|^2 \leq \liminf_{m \to \infty} R_m^2 - |v - \overline{u}_m|^2 \]
so $(\overline{u}_{\infty}, R_{\infty}) \in \mathcal{B}$ realizes the infimum. 

Assume now by contradiction that $R_m \to \infty$, and so $|\overline{u}_m| \to \infty$ as well. Notice that 
\[ R_m^2 - |v - \overline{u}_m|^2 = (R_m - |v - \overline{u}_m|) (R_m + |v - \overline{u}_m|) \]
But $R_m + |v - \overline{u}_m| \to \infty$, so to get convergence it is needed that $R_m - |v - \overline{u}_m| \to 0$. However $v \notin \partial \mathcal{K}$ so there exists $\delta > 0$ such that $B(v, \delta) \subset \mathcal{K} \subset B(\overline{u}_m, R_m)$ for every $m$ (the second inclusion coming from the fact that $(\overline{u}_m, R_m) \in \mathcal{B}$), therefore $|v - \overline{u}_m| \leq R_m - \delta$, which is a contradiction. 

So the infimum has to be realized. In particular, 
\begin{equation} \inf_{(\overline{u}, R) \in \mathcal{B}} \frac{|v|^2 + R^2 - |\overline{u}-v|^2}{n} ~ > \lambda_1 \label{equ-preuvestep2-lambda1} \end{equation}
Let us then consider $v_t = v + t e_1$, and $M_t = M + t v \otimes e_1 + t e_1 \otimes v - 2 t (v \cdot e_1) e_1 \otimes e_1$. It is easy to check that $M_t$ is trace-free for every $t$. Furthermore,  
\[ v_t \otimes v_t - M_t = v \otimes v - M + t^2 e_1 \otimes e_1 - 2 t (v \cdot e_1) e_1 \otimes e_1 \]
In particular, the eigenvalues of $v_t \otimes v_t - M_t$ are the same as those of $v \otimes v - M$, except maybe $\lambda_1$. But since \eqref{equ-preuvestep2-lambda1} is open, we deduce that \eqref{equ-preuvestep2-lambda1} remains true for small (positive or negative) values of $t$. On the other hand, for small values of $t$, $v_t$ remains in $\mathcal{K}$ because $v$ is in the interior of $\mathcal{K}$. So by setting $q_t = q$, we have for any small enough $t$ that $(v_t, M_t, q_t) \in C_{\mathcal{K}}$ and it is linear in $t$, non-constant. This means that $(v, M, q)$ is not an extremal point. 
\end{itemize} 

\paragraph{Step 3 : Identification of the convex hull.} By the previous steps, $C_{\mathcal{K}}$ is a convex set, that contains $\overline{\mathcal{K}}$ and $\overline{\mathcal{K}}$ contains all extremal points of $C_{\mathcal{K}}$. This means $C_{\mathcal{K}} = \overline{\mathcal{K}}^{co}$ by the Krein-Milman theorem. 

Moreover, we proved that all extremal points satisfy $v \in \partial \mathcal{K}$. In particular, this implies $\overline{\mathcal{K}}^{co} \cap \left( \partial \mathcal{K} \times \mathcal{S}^n_0 \times \mathbb{R} \right) \subset \overline{\mathcal{K}}$, and the converse is automatic. 

\paragraph{Step 4 : Identification of the $\Lambda$-convex hull.} Let now $(v, M, q) \in C_{\mathcal{K}} = \overline{\mathcal{K}}^{co}$. We can then write : 
\[ (v, M, q) = \sum_{i = 1}^N \alpha_i (v_i, M_i, q_i) \]
for some points $(v_i, M_i, q_i) \in \overline{\mathcal{K}}$, by definition of the convex hull, with $(\alpha_i)$ positive coefficients satisfying $\sum_{i = 1}^N \alpha_i = 1$ and $N \geq 1$. Note that, if $(v_i, M_i, q_i) \in \overline{\mathcal{K}}$, then $(v_i, M_i, q_i') \in \overline{\mathcal{K}}$ for every $q_i' \in \mathbb{R}$. We may thus choose $q_i = q$ for every $i$. Furthermore, belonging to $\overline{\mathcal{K}}$ forces $M_i = F(v_i)$. If $N = 1$, we have $(v, M, q) \in \overline{\mathcal{K}}$ which lies indeed in the $\Lambda$-convex hull of $\overline{\mathcal{K}}$. Assume now $N \geq 2$. 

We procede by induction and assume that, for every $(v', M', q') \in \overline{\mathcal{K}}^{co}$, if $(v', M', q')$ can be written as barycenter of at most $(N-1)$ points of $\overline{\mathcal{K}}$, then $(v', M', q')$ belongs to the $\Lambda$-convex hull of $\overline{\mathcal{K}}$. 

By lemma \ref{lem-lambda-difference}, there exists $q_{12} \in \mathbb{R}$ such that $(v_2 - v_1, F(v_2) - F(v_1), q_{12}) \in \Lambda$. We then write : 
\[ \begin{aligned}
(v, M, q) &= \alpha_1 (v_1, F(v_1), q) + \alpha_2 (v_2, F(v_2), q) + \sum_{i = 3}^N \alpha_i (v_i, F(v_i), q) \\
&= \frac{\alpha_1}{\alpha_1 + \alpha_2} \left( (\alpha_1 + \alpha_2) \left( v_1, F(v_1), q - \frac{q_{12}}{2} \right) + \sum_{i = 3}^N \alpha_i (v_i, F(v_i), q) \right) \\
&\quad \quad \quad + \frac{\alpha_2}{\alpha_1 + \alpha_2} \left( (\alpha_1 + \alpha_2) \left( v_2, F(v_2), q + \frac{q_{12}}{2} \right) + \sum_{i = 3}^N \alpha_i (v_i, F(v_i), q) \right) \\
&=: \frac{\alpha_1}{\alpha_1 + \alpha_2} (v', M', q') + \frac{\alpha_2}{\alpha_1 + \alpha_2} (v'', M'', q'') 
\end{aligned} \]
Notice that $(v', M', q')$ and $(v'', M'', q'')$ are two points of $\overline{\mathcal{K}}^{co}$ that can be written as barycenter of $N-1$ points of $\overline{\mathcal{K}}$, so they both lie in the $\Lambda$-convex hull of $\overline{\mathcal{K}}$ by induction hypothesis. But $(v'', M'', q'') - (v', M', q') = (v_2 - v_1, F(v_2) - F(v_1), q_{12}) \in \Lambda$, so $(v, M, q)$ also belongs to the $\Lambda$-convex hull of $\overline{\mathcal{K}}$. 

By induction, we deduce that $\overline{\mathcal{K}}^{co}$ is also the $\Lambda$-convex hull of $\overline{\mathcal{K}}$. 

\paragraph{Step 5 : Good oscillation directions.} We now prove the last statement of the proposition. To that end, let $z \in \mbox{int} ~ \overline{\mathcal{K}}^{co}$. As in step 4, we then write : 
\[ z = (v, M, q) = \sum_{i = 1}^N \alpha_i (v_i, F(v_i), q) \]
The $(\alpha_i)$ are chosen in $(0, 1)$ and satisfy $\sum_{i = 1}^N \alpha_i = 1$. 

\begin{itemize} \item First, let us prove that we may choose the $(v_i)$ such that, for every $(i, j)$, we have $v_i$ and $v_j$ not aligned (that is, $|v_i|^2 |v_j|^2 - (v_i \cdot v_j)^2 \neq 0$). Indeed, if given $(e_i)_{1 \leq i \leq D}$ an orthonormal basis of $\mathbb{R}^n \times \mathcal{S}^n_0 \times \mathbb{R}$ and $\delta > 0$ small enough such that $B(z, D \delta) \subset \overline{\mathcal{K}}^{co}$, then if we set $z_i = z + \delta e_i \in \mbox{int} ~ \overline{\mathcal{K}}^{co}$, and $z_0 = z - \delta \sum_i e_i \in \mbox{int} ~ \overline{\mathcal{K}}^{co}$, we may write for every $0 \leq i \leq D$, 
\[ z_i = \sum_{j = 1}^{N_i} \beta_{ij} (v_{ij}, F(v_{ij}), q_{ij}) \]
with positive $(\beta_{ij})$, summing to $1$, and $N_i \leq D + 1 = N_{max}(n)$ by Carathéodory's theorem. We may also choose $v_{ij}$ to be extremal points, so $v_{ij} \in \partial \mathcal{K}$. 

However, the map
\[ (v_{ij})_{0 \leq i \leq D, 1 \leq j \leq N_i} \in (\partial \mathcal{K})^{N_0 + ... + N_D} \mapsto \left( \sum_{j = 1}^{N_i} \beta_{ij} (v_{ij}, F(v_{ij}), q_{ij}) \right)_{0 \leq i \leq D} \]
is continuous and, since $\mathcal{K}$ is a convex compact set of non-empty interior, $\partial \mathcal{K}$ is homeomorphic to $S^{n-1}$. It is clear that $z$ belongs to the interior of the convex hull of the $(z_i)$, so we may perturb the $(v_{ij})$ on $\partial \mathcal{K}$ in order to satisfy $v_{ij} + v_{i'j'} \neq 0$ for every $i, i', j, j'$, and so obtain the desired expression of $z$. By Caratheodory's theorem, we end up with an expression of $z$ as a convex combination of $N'$ points, with $N'$ being bounded by a universal constant depending only on the dimension $n$. 

\item Assume now $\alpha_1 \geq \alpha_2 \geq ... \geq \alpha_N > 0$. Let then $j > 1$ be such that
\[ \alpha_j |v_j - v_1| = \max_{i > 2} \alpha_i |v_i - v_1| \]
Then we set
\[ \lambda = \frac{\alpha_j}{2} \left( v_j - v_1, F(v_j) - F(v_1), \frac{|v_j|^2 - |v_1|^2}{n} \right) = (\overline{v}, \overline{M}, \overline{q}) \]
For every $t \in [-1, 1]$, 
\[ (v, M, q) + t\lambda = \sum_{i = 1}^N \alpha_i(t) (v_i, F(v_i), q_i(t)) \]
for $\alpha_i(t) = \alpha_i + \frac{t \alpha_j}{2} (\delta_{i, j} - \delta_{i, 1})$, and $q_i(t) = q$ for $i \neq 1, j$. Note that $\alpha_i(t) \in [0, 1]$ for any $t \in [-2, 2]$, and $\sum_{i = 1}^N \alpha_i(t) = 1$ for any $t$. Indeed, if $i \neq 1$ and $i \neq j$, $\alpha_i(t) = \alpha_i \in (0, 1)$ for any $t$ ; for $i = 1$, 
\[ \alpha_1(t) = \alpha_1 - \frac{t \alpha_j}{2} \geq 0 \]
if $t \leq -2$ since $\alpha_1 \geq \alpha_j > 0$, and so $\alpha_j(t) \leq 1 - \alpha_1(t) \leq 1$ ; for $i = j$, 
\[ \alpha_j(t) = \alpha_j + \frac{t \alpha_j}{2} \geq 0 \]
if $t \geq 2$, and so $\alpha_1(t) \leq 1 - \alpha_j(t) \leq 1$. 

Therefore, $(v, M, q) + t\lambda \in \overline{\mathcal{K}}^{co}$ for every $t \in [-2, 2]$. Moreover, 
\[ \mbox{dist}(v, \partial \mathcal{K}) \leq |v - v_1| = \left| \sum_{i = 1}^N \alpha_i (v_i - v_1) \right| \leq (N-1) \alpha_j |v_j - v_1| = 2 (N-1) |\overline{v}| \]
We already explained why $N$ was bounded by a function of $n$, so we get the second inequality of \eqref{equ-contdist-propenvconv}. 

The other inequality comes from convexity. Indeed, let $t \in [0, 1]$. By definition, if $r = \mbox{dist}((v, M, q), \partial(\overline{\mathcal{K}}^{co}))$, then $B((v, M, q), r) \subset \overline{\mathcal{K}}^{co}$. The choice of $\lambda$ is such that $(v, M, q) + 2 \lambda \in \overline{\mathcal{K}}^{co}$. By convexity, this means that 
\[ B\left( (v, M, q) + t \lambda, \frac{2-t}{2}r \right) \subset \overline{\mathcal{K}}^{co} \]
which implies that 
\begin{align*}
\mbox{dist}((v, M, q) + t \lambda, \partial(\overline{\mathcal{K}}^{co})) \geq \frac{2-t}{2} r ~ \geq \frac{1}{2} \mbox{dist}((v, M, q), \partial(\overline{\mathcal{K}}^{co})) 
\end{align*}
We can procede likewise for $t \in [-1, 0]$. 
\end{itemize} 
\end{Dem}

We may use essentially the same proof in the case of translated sets $\overline{\mathcal{K}(z_0)}$ : 

\begin{Cor} Let $\mathcal{K} \in \mathfrak{K}$ be uniformly convex, $z_0 = (v_0, M_0, q_0) \in \mathbb{R}^n \times \mathcal{S}^n_0(\mathbb{R}) \times \mathbb{R}$ and $r > 0$. Then proposition \ref{prop-enveloppe-convexe} also applies with $\overline{\mathcal{K}(z_0, r)}$ : the $\Lambda$-convex hull of $\overline{\mathcal{K}(z_0, r)}$ is the same as the convex hull of $\overline{\mathcal{K}(z_0, r)}$, and it is equal to
\[ \begin{aligned} \overline{\mathcal{K}(z_0, r)}^{co} =& \left\{ (v, M, q) \in (v_0 + r\mathcal{K}) \times \mathcal{S}^n_0 \times \mathbb{R}, ~~ \forall (\overline{u}, R) \in \mathbb{R}^n \times \mathbb{R}^{+*} \right. \\
&\left. \mbox{ such that } r\mathcal{K} \subset B(\overline{u}, R), ~~ v \otimes v - M \leq \frac{|v|^2 + R^2 - |v - v_0 - \overline{u}|^2}{n} I_n \right\} 
\end{aligned} \]
Furthermore, $\overline{\mathcal{K}(z_0, r)} = \overline{\mathcal{K}(z_0, r)}^{co} \cap \{ (v, M, q), v - v_0 \in r \partial \mathcal{K} \}$. 

Moreover, with the same constant $C = C(n) > 0$, for every $z = (v, M, q)$ in the interior of $\overline{\mathcal{K}(z_0, r)}^{co}$, there exists $\lambda = (\overline{v}, \overline{M}, \overline{q}) \in \Lambda_{\mathcal{K}}(z_0, r)$ such that $[z-\lambda, z+\lambda] \in \mbox{int} ~ \overline{\mathcal{K}(z_0, r)}^{co}$ and 
\[ \mbox{dist}\left( \left[ z - \lambda, z + \lambda \right], \partial (\overline{\mathcal{K}(z_0, r)}^{co}) \right) \geq \frac{1}{2} d(z, \partial (\overline{\mathcal{K}(z_0, r)}^{co})), \quad |\overline{v}| \geq C \mbox{dist}(v-v_0, r\partial \mathcal{K}) \] \label{cor-enveloppe-convexe}
\end{Cor}

\begin{Dem}
Let us first notice that
\begin{align*} 
\overline{\mathcal{K}(z_0, r)} &= \left\{ (v, M, q) \in K, ~ v-v_0 \in r\partial \mathcal{K} \right\} \\
&= \left\{ (v_0 + \widetilde{v}, F(v_0) + v_0 \otimes \widetilde{v} + \widetilde{v} \otimes v_0 - \frac{2 v_0 \cdot \widetilde{v}}{n} I_n + F(\widetilde{v}), q), ~ \widetilde{v} \in r \partial \mathcal{K}, q \in \mathbb{R} \right\} \end{align*}
In particular, if we consider the affine and invertible map 
\[ G : (v, M, q) \mapsto \left( r v + v_0, r^2 M + F(v_0) + v_0 \otimes (rv) + (rv) \otimes v_0 - \frac{2 v_0 \cdot (rv)}{n} I_n, r^2 q + \frac{|v_0|^2}{n} + \frac{2 v_0 \cdot (rv)}{n} \right) \]
we have that
\[ \overline{\mathcal{K}(z_0, r)} = G(\overline{\mathcal{K}}) \]
so 
\[ \begin{aligned}
&\overline{\mathcal{K}(z_0, r)}^{co} = G(\overline{\mathcal{K}}^{co}) \\
&= \Big\{ (v, M, q) \in (v_0 + r\mathcal{K}) \times \mathcal{S}^n_0 \times \mathbb{R}, \forall (\overline{u}, R) \in \mathbb{R}^n \times \mathbb{R}^{+*} \mbox{ such that } r\mathcal{K} \subset B(\overline{u}, R), \\
&(v - v_0) \otimes (v - v_0) - \left( M - F(v_0) - v_0 \otimes (v - v_0) - (v - v_0) \otimes v_0 + \frac{2 v_0 \cdot (v- v_0)}{n} I_n \right) \\
&\quad \quad \quad \quad \quad \leq \frac{|v-v_0|^2 + R^2 - |v-v_0 - \overline{u}|^2}{n} I_n \Bigl\} \\
&= \left\{ (v, M, q) \in (v_0 + r\mathcal{K}) \times \mathcal{S}^n_0 \times \mathbb{R}, \forall (\overline{u}, R) \in \mathbb{R}^n \times \mathbb{R}^{+*} \right. \\
&\left. \mbox{ such that } r\mathcal{K} \subset B(\overline{u}, R), v \otimes v - M \leq \frac{|v|^2 + R^2 - |v-v_0 - \overline{u}|^2}{n} I_n \right\} 
\end{aligned} \]

This replaces steps 1 to 3. Step 4 can be done the same way. Finally, for step 5, we notice that the linear part of $G$ maps $\Lambda_{\mathcal{K}}$ in $\Lambda_{\mathcal{K}}(z_0, r)$ up to the non-alignment condition ($a$ and $b$ always remain distinct, but can become aligned after a translation by $v_0$). However, using the same perturbative argument as in step 5, we may avoid this as well and conclude. 
\end{Dem}

\section{Localized oscillating solutions} \label{section-localized}

\subsection{Potential}

Let $\mathcal{K} \in \mathcal{K}$, $z_0$ a base point, $r > 0$. 

\begin{Prop} For every $\lambda \in \Lambda_{\mathcal{K}}(z_0, r)$, there exists an order $3$ homogeneous differential operator with constant coefficients, 
\[ A(\nabla) : C^{\infty}_c(\mathbb{R}^{n+1}) \to C^{\infty}_c(\mathbb{R}^{n+1}, \mathbb{R}^n \times \mathcal{S}^n_0 \times \mathbb{R}) \]
and $(c, \xi) \in \mathbb{R} \times \mathbb{R}^n$, $\xi \neq 0$, such that
\begin{enumerate}
\item for every $\varphi \in C^{\infty}_c(\mathbb{R}^{n+1})$, $A(\nabla) \varphi$ is a solution of the differential system \eqref{equ-Euler-lin}
\item if $\psi \in C^{\infty}_c(\mathbb{R})$ and that we set $\varphi(t, x) = \psi\left( ct + \xi \cdot x \right)$, then 
\[ \left[ A(\nabla) \varphi \right] (t, x) = \lambda \psi'''\left( ct + \xi \cdot x \right) \]
\end{enumerate} \label{prop-solutions-potentiel}
\end{Prop}

\begin{Dem}
Let $\lambda \in \Lambda_{\mathcal{K}}(z_0, r)$. Let then $t \in \mathbb{R}^{+}$ and $a, b \in v_0 + r\partial \mathcal{K}$ be such that $\lambda = t(v, M, q)$, $v = b-a$, $M = F(b) - F(a)$, $q = \frac{|b|^2 - |a|^2}{n}$, and $(a, b)$ is a linearly independent family (cf definition \ref{def-Lambda-de-K}). 

Denote $A(\nabla) = (A_v(\nabla), A_M(\nabla), A_q(\nabla))$ and let us set 
\[ \left\{ \begin{array}{l}
A_v^i = \sum_{k, j} (b^k b^j - a^k a^j) \partial_{kji} + \sum_{k, j} (a^k a^i - b^k b^i) \partial_{kjj}, \\
A_M^{ij} = \sum_k (b^i b^j - a^i a^j) \partial_{tkk} + (|a|^2 - |b|^2) \frac{1}{n} \delta^{ij} \partial_t \Delta_x, \\
A_q = (|b|^2 - |a|^2) \frac{1}{n} \partial_t \Delta_x + \sum_{k, j} (a^k a^j - b^k b^j) \partial_{tkj} 
\end{array} \right. \]
where $\delta^{ij}$ is the Kronecker symbol, and $\partial_j = \partial_{x_j}$. In vector notation, 
\[ \left\{ \begin{array}{l}
A_v = \nabla_x (b \cdot \nabla_x) (b \cdot \nabla_x) - \nabla_x (a \cdot \nabla_x) (a \cdot \nabla_x) + a (a \cdot \nabla_x) \Delta_x - b (b \cdot \nabla_x) \Delta_x \\
A_M = (b \otimes b) \partial_t \Delta_x - (a \otimes a) \partial_t \Delta_x + (|a|^2 - |b|^2) \frac{1}{n} I_n \partial_t \Delta_x \\
A_q = (|b|^2 - |a|^2) \frac{1}{n} \partial_t \Delta_x + \partial_t (a \cdot \nabla_x) (a \cdot \nabla_x) - \partial_t (b \cdot \nabla_x) (b \cdot \nabla_x) 
\end{array} \right. \]

We then check the first statement : 
\[ \nabla_x \cdot A_v(\nabla) = \Delta_x (b \cdot \nabla_x)^2 - \Delta_x (a \cdot \nabla_x)^2 + (a \cdot \nabla_x)^2 \Delta_x - (b \cdot \nabla_x)^2 \Delta_x = 0 \]
and 
\[ \begin{aligned}
&\partial_t A_v(\nabla) + \nabla_x \cdot A_M(\nabla) + \nabla_x A_q(\nabla) \\
&= \partial_t \nabla_x (b \cdot \nabla_x)^2 - \partial_x \nabla_x (a \cdot \nabla_x)^2 + \partial_t a (a \cdot \nabla_x) \Delta_x - \partial_t b (b \cdot \nabla_x) \Delta_x \\
&+ b (b \cdot \nabla_x) \partial_t \Delta_x - a (a \cdot \nabla_x) \partial_t \Delta_x + (|a|^2 - |b|^2) \frac{1}{n} \nabla_x\partial_t \Delta_x \\
&+ (|b|^2 - |a|^2) \frac{1}{n} \nabla_x \partial_t \Delta_x + \partial_t \nabla_x (a \cdot \nabla_x)^2 - \partial_t \nabla_x (b \cdot \nabla_x)^2 \\
&= 0
\end{aligned} \]
Furthermore, 
\[ \mbox{Tr}(A_M(\nabla)) = |b|^2 \partial_t \Delta_x - |a|^2 \partial_t \Delta_x + (|a|^2 - |b|^2) \partial_t \Delta_x = 0 \]
and $A_M(\nabla)^T = A_M(\nabla)$ so $A_M(\nabla)$ indeed maps $C_c^{\infty}(\mathbb{R}^{n+1}, \mathbb{R})$ into $C_c^{\infty}(\mathbb{R}^{n+1}, \mathcal{S}^n_0(\mathbb{R}))$. 

Then, let us set $\delta \in \mathbb{R}$ the only (real) solution of 
\[ \delta^3 = -\frac{t}{(|b|^2 |a|^2 - (a \cdot b)^2)^2 |b-a|^2} \]
which is well defined if $(a, b)$ is a linearly independent family. We may assume $t > 0$ : otherwise, $\lambda = 0$ and the lemma is trivial by setting $A(\nabla) = 0$. Therefore, $\delta \neq 0$ and we set 
\[ \xi = \delta \left( (|b|^2 - a \cdot b) a + (|a|^2 - a \cdot b) b \right), \quad c = -\xi \cdot a \]
A computation shows that
\[ c = -\xi \cdot a = - \delta \left( (|b|^2 - a \cdot b) |a|^2 + (|a|^2 - a \cdot b) a \cdot b \right) = - \delta \left( |b|^2|a|^2 - (a \cdot b)^2 \right) \neq 0 \]
and so $\xi \neq 0$. The symmetry of this last expression of $c$ in $a, b$ (as well as of the expression of $\xi$) implies that $c = - \xi \cdot b$. We may also compute 
\[ \begin{aligned}
|\xi|^2 &= \delta^2 \left( \left( |b|^2 - a \cdot b \right)^2 |a|^2 + \left( |a|^2 - a \cdot b \right)^2 |b|^2 + 2 (a \cdot b) \left( |b|^2 - a \cdot b \right) \left( |a|^2 - a \cdot b \right) \right) \\
&= \delta^2 \Big( |b|^4 |a|^2 + (a \cdot b)^2 |a|^2 - 2 |a|^2 |b|^2 (a \cdot b) + |a|^4 |b|^2 + (a \cdot b)^2 |b|^2 - 2 |a|^2 |b|^2 (a \cdot b) \\
&\quad \quad \quad \quad + 2 (a \cdot b) |a|^2 |b|^2 + 2 (a \cdot b)^3 - 2 |b|^2 (a \cdot b)^2 - 2 |a|^2 (a \cdot b)^2 \Big) \\
&= \delta^2 \left( |b|^4 |a|^2 + |a|^4 |b|^2 - (a \cdot b)^2 |a|^2 - 2 |a|^2 |b|^2 (a \cdot b) - (a \cdot b)^2 |b|^2 + 2 (a \cdot b)^3 \right) \\
&= \delta^2 \left( |b|^2 |a|^2 (|a|^2 + |b|^2 - 2 (a \cdot b)) - (a \cdot b)^2 (|a|^2 + |b|^2 - 2 (a \cdot b)) \right) \\
&= \delta^2 |a - b|^2 \left( |b|^2 |a|^2 - (a \cdot b)^2 \right) = - \delta |a - b|^2 c 
\end{aligned} \]
Finally, we have the relation
\[ c^2 \delta |b-a|^2 = -t \]

Let then $\psi \in C^{\infty}_c(\mathbb{R})$ and let us set $\varphi(t, x) = \psi\left( c t + \xi \cdot x \right)$. Then we can compute 
\[ \begin{aligned}
A_v(\nabla) \varphi(t, x) &= \left( \xi (b \cdot \xi)^2 - \xi (a \cdot \xi)^2 + a (a \cdot \xi) |\xi|^2 - b (b \cdot \xi) |\xi|^2 \right) \psi'''\left( ct + \xi \cdot x \right) \\
A_M(\nabla) \varphi(t, x) &= c |\xi|^2 \left( (b \otimes b) - (a \otimes a) + (|a|^2 - |b|^2) \frac{1}{n} I_n \right) \psi'''\left( ct + \xi \cdot x \right) \\
A_q(\nabla) \varphi(t, x) &= \left( (|b|^2 - |a|^2) \frac{1}{n} c |\xi|^2 + c (a \cdot \xi)^2 - c (b \cdot \xi)^2 \right) \psi'''\left( ct + \xi \cdot x \right) 
\end{aligned} \]
Yet we may also compute
\[ \begin{aligned}
&\xi (b \cdot \xi)^2 - \xi (a \cdot \xi)^2 + a (a \cdot \xi) |\xi|^2 - b (b \cdot \xi) |\xi|^2 \\
&\quad = c^2 \delta \left( \xi - \xi + a |a - b|^2 - b |a - b|^2 \right) \\
&\quad = c^2 \delta (a - b) |a - b|^2 = t(b - a) \\
&c |\xi|^2 \left( (b \otimes b) - (a \otimes a) + (|a|^2 - |b|^2) \frac{1}{n} I_n \right) \\
&\quad = - c^2 \delta |a - b|^2 \left( F(b) - F(a) \right) = t \left( F(b) - F(a) \right) \\
&(|b|^2 - |a|^2) \frac{1}{n} c |\xi|^2 + c (a \cdot \xi)^2 - c (b \cdot \xi)^2 \\
&\quad = c^2 \delta |a - b|^2 (|a|^2 - |b|^2) \frac{1}{n} + c^3 - c^3 = t \frac{|b|^2 - |a|^2}{n} 
\end{aligned} \]
So we deduce that
\[ A(\nabla) \varphi(t, x) = \lambda \psi'''\left( ct + \xi \cdot x \right) \]
as desired. 
\end{Dem}

\subsection{Localisation}

\begin{Prop} Let $O$ be a bounded non empty open set of $\mathbb{R}^n$, $I = (t_0, t_1) \subset \mathbb{R}^n$, $\mathcal{K} \in \mathfrak{K}$, $z_0$ a base point, $r > 0$ $\lambda = (\overline{v}, \overline{M}, \overline{q}) \in \Lambda_{\mathcal{K}}(z_0, r)$, and $\mathcal{V}$ a neighborhood of $[-\lambda, \lambda] \subset \mathbb{R}^n \times \mathcal{S}^n_0 \times \mathbb{R}$. Let $O'$ be a non empty open set $\mathbb{R}^n$ such that $\overline{O'} \subset O$, let $\theta \in \left( 0, \frac{t_1 - t_0}{2} \right)$ and set $I_{\theta} = [t_0 + \theta, t_1 - \theta]$. Let also $\mathbf{K} : \overline{I} \times \overline{O} \to \mathfrak{K}$ be continuous for the Hausdorff distance, $u : \overline{I} \times \overline{O} \to \mathbb{R}^n$ be continuous, and assume that there exists $c_0 > 0, R_0 > r_0 > 0$ such that for every $(t, x) \in \overline{I} \times \overline{O}$, $\mathbf{K}(t, x)$ is $c_0$-uniformly convex, $r_{min}(\mathbf{K}(t, x)) \geq r_0$, $r_{max}(\mathbf{K}(t, x)) \leq R_0$. 

Then for any $\eta > 0$, there exists $z = (v, M, q) \in C^{\infty}_c(I \times O, \mathcal{V})$ such that
\begin{enumerate}
\item $(v, M, q)$ is a solution of the differential system \eqref{equ-Euler-lin}
\item for every $t \in I_{\theta}$, 
\[ \int_O j_{\mathbf{K}(t, x)}(u(t, x) + v(t, x))^2 ~ dx ~~ \geq \int_O j_{\mathbf{K}(t, x)}(u(t, x))^2 ~ dx ~ + \frac{c_0 r_0^2}{6} |O'| |\overline{v}|^2 \]
\item $\sup_{t \in I} \Vert z(t) \Vert_{H^{-1}(O)} < \eta$
\end{enumerate} \label{prop-solutions-loc}
\end{Prop}

\begin{Dem}
The beginning of the proof is the same as in \cite{VillaniBourbaki}, see also \cite{DeLellisSzekelyhidi2}. 

Given such a $\lambda$, let $A, \xi, c$ be as in proposition \ref{prop-solutions-potentiel} : in particular, $\xi \neq 0$. Let $\varphi \in C^{\infty}_c(\mathbb{R}^{n+1}, [0, 1])$ be a function supported in $I \times O$ and identically $1$ on $I_{\theta} \times O'$. Then, for any $\varepsilon > 0$, consider the functions 
\[ \varphi_{\varepsilon}(t, x) := \varepsilon^3 \varphi(t, x) \cos\left( \frac{c t + \xi \cdot x}{\varepsilon} \right) \]
and
\[ z_{\varepsilon}(t, x) = A(\nabla) \varphi_{\varepsilon}(t, x) \]
Let us also define
\[ \phi_{\varepsilon}(t, x) = \varepsilon^3 \cos\left( \frac{ct + \xi \cdot x}{\varepsilon} \right), \quad \psi_{\varepsilon}(s) = \varepsilon^3 \cos\left( \frac{s}{\varepsilon} \right) \]
By proposition \ref{prop-solutions-potentiel}, $z_{\varepsilon}$ is a solution of \eqref{equ-Euler-lin}, belonging to $C^{\infty}_c$ and compactly supported in $I \times O$ (as $\varphi_{\varepsilon}$). 

Moreover, 
\[ z_{\varepsilon}(t, x) = \varphi(t, x) A(\nabla) \phi_{\varepsilon}(t, x) + O(\varepsilon) \]
when $\varepsilon \to 0$, where the $O$ is uniform in $t, x$, and by proposition \ref{prop-solutions-potentiel}, $A(\nabla) \phi_{\varepsilon}(t, x) = \lambda \psi_{\varepsilon}'''\left( ct + x \cdot \xi \right)$. But $\Vert \psi_{\varepsilon}''' \Vert_{L^{\infty}} = 1$ so by choosing $\varepsilon$ small enough we get that $z_{\varepsilon}$ is valued in $\mathcal{V}$. 

Let us now prove the second point. Note that for any fixed $t \in I_{\theta}$, since 
\[ v_{\varepsilon}(t, x) = \varphi(t, x) \overline{v} \sin\left( \frac{c t + \xi \cdot x}{\varepsilon} \right) + O(\varepsilon) \]
and $\xi \neq 0$, one has by continuity of $\mathbf{K}, u, \varphi$ and corollary \ref{lem-jauge-Lipschitz-tot}, that
\begin{equation}
\begin{aligned} 
\int_O j_{\mathbf{K}(t, x)}(u(t, x) + v_{\varepsilon}(t, x))^2 ~ dx ~ \underset{\varepsilon \to 0}{\longrightarrow} &~ \fint_0^{2 \pi} \int_O j_{\mathbf{K}(t, x)}(u(t, x) + \varphi(t, x) \overline{v} \sin(y))^2 ~ dx dy \\
&\geq \int_O j_{\mathbf{K}(t, x)}(u(t, x))^2 ~ dx ~ \\
&\quad \quad \quad + ~ \frac{c_0}{2} \int_O \fint_0^{2\pi} |\varphi(t, x)|^2 |\sin(y)|^2 j_{\mathbf{K}(t, x)}(\overline{v})^2 ~ dy dx \\
&\geq \int_O j_{\mathbf{K}(t, x)}(u(t, x))^2 ~ dx ~ + ~ \frac{c_0}{4} \int_{O'} j_{\mathbf{K}(t, x)}(\overline{v})^2 ~ dx \\
&\geq \int_O j_{\mathbf{K}(t, x)}(u(t, x))^2 ~ dx ~ + ~ \frac{c_0 r_0^2}{4} |O'| |\overline{v}|^2
\end{aligned} \label{equ-lemloc-min}
\end{equation}
by $c_0$-uniform convexity of $\mathbf{K}$. 

Let us show that this convergence is uniform in $t \in I_{\theta}$. It is clear that the right-hand side is continuous in $t$ under our assumptions. On the left-hand side, notice that $v_{\varepsilon}$ can be written as a finite sum of functions of the form 
\[ f_1(t, x) f_2\left( \frac{ct+x\cdot\xi}{\varepsilon} \right) \]
for $f_1$ a function independant of $\varepsilon$ localizing near $I_{\theta} \times O'$ and $f_2 \in \{ \cos, \sin \}$. In particular, if $|\eta| < \frac{\theta}{2}$, 
\[ f_1(t+\eta, x) f_2\left( \frac{c(t+\eta) + x \cdot \xi}{\varepsilon} \right) = f_1(t+\eta, x') f_2\left( \frac{ct + x' \cdot \xi}{\varepsilon} \right) \]
for $x' = x + \alpha \xi$ for $\alpha \in [-\varepsilon \pi, \varepsilon \pi)$ unique such that $\alpha |\xi|^2 = c\eta$. 

Therefore, applying a change of variable $x \mapsto x'$ ($\varepsilon$ being fixed), we obtain that 
\begin{align*}
&\int_O j_{\mathbf{K}(t+\eta, x)}(u(t+\eta, x)+v_{\varepsilon}(t+\eta, x))^2 ~ dx \\
&= \int_{O+\alpha \xi} j_{\mathbf{K}(t+\eta, x-\alpha \xi)}\left(u(t+\eta, x - \alpha \xi) + \sum f_1(t+\eta, x-\alpha \xi) f_2\left( \frac{ct + x \cdot \xi}{\varepsilon} \right) \right)^2 ~ dx
\end{align*}
But by lemma \ref{lem-jauge-Lipschitz-tot}, we have that 
\begin{align*} 
j_{\mathbf{K}(t+\eta, x-\alpha \xi)}&\left(u(t+\eta, x - \alpha \xi) + \sum f_1(t+\eta, x-\alpha \xi) f_2\left( \frac{ct + x \cdot \xi}{\varepsilon} \right) \right)^2 \\
&\underset{\eta \to 0}{\longrightarrow} j_{\mathbf{K}(t, x)}\left(u(t, x) + \sum f_1(t, x) f_2\left( \frac{ct + x \cdot \xi}{\varepsilon} \right) \right)^2 
\end{align*}
for every $(x, t)$ (note that $\alpha = \alpha(\eta)$ goes to $0$ as $\eta \to 0$) ; moreover, this convergence is uniform in $\varepsilon$ since $|\alpha| \leq C|\eta|$ for some $C > 0$, uniformly in $\varepsilon$. We deduce that 
\begin{align*} \int_O j_{\mathbf{K}(t+\eta, x)}(u(t+\eta, x)+v_{\varepsilon}(t+\eta, x))^2 ~ dx \underset{\eta \to 0}{\longrightarrow} \int_O j_{\mathbf{K}(t, x)}(u(t, x)+v_{\varepsilon}(t, x))^2 ~ dx \end{align*}
uniformly in $\varepsilon > 0$ by compactness of $\overline{O}$ and continuity of all the functions over $\overline{O}$. 

Therefore, since both the left-hand side and the right-hand side of \eqref{equ-lemloc-min} are continuous in $t \in I_{\theta}$ compact, uniformly in $\varepsilon > 0$, we deduce that that for $\varepsilon > 0$ small enough the second point of the lemma holds. 

Finally, for the third and last point, we have for any $\zeta \in H^1$ with $\Vert \zeta \Vert_{H^1} \leq 1$ and any $t \in I$ that after one integration by parts : 
\begin{align*}
\int_O z(t, x) \zeta(x) ~ dx = O(\varepsilon) \int_O |\nabla \zeta(x)| ~ dx \leq O(\varepsilon) |O| \Vert \zeta \Vert_{H^1}
\end{align*}
where the $O(\varepsilon)$ is uniform in $t$. By choosing $\varepsilon > 0$ small enough with respect to $\eta$, $O$ being bounded, we conclude. 
\end{Dem}

\section{Convex integration} \label{section-convexint}

In this section, we prove theorem \ref{theo-principal}. We fix $z_0 = (v_0, M_0, q_0)$, $\mathbf{K}$, $c_0, r_0, R_0$ as in the theorem. 

\subsection{The space of subsolutions}

Set 
\[ \begin{aligned}
X_0 &= \Biggl\{ z = (v, M, q) \in \left( v_0, M_0, q_0 \right) + C^{\infty}_c(\mathbb{R} \times \mathbb{R}^n, \mathbb{R}^n \times \mathcal{S}^n_0 \times \mathbb{R}), ~ \mbox{ solution of } \eqref{equ-Euler-lin}, \\
&\forall (t, x) \in \mathbb{R}^{n+1}, z(t, x) \in ~ \mbox{int} ~ \overline{\mathbf{K}(t, x)(z_0(t, x), \sqrt{\overline{e}(t, x)})}^{co} \mbox{ or } (v(t, x), M(t, x)) = (v_0(t, x), M_0(t, x)) \Biggl\} 
\end{aligned} \]
This is the set of subsolutions. 

Notice that by corollary \ref{cor-enveloppe-convexe} and the hypothesis of theorem \ref{theo-principal}, we have for every $(t, x) \in \mathbb{R}^{n+1}$ not satisfying $F(v_0(t, x)) = M_0(t, x)$ that
\[ F(v_0(t, x)) - M_0(t, x) < \overline{e}(t, x) \frac{a_{\mathbf{K}(t, x)}}{n} I_n = \frac{a_{\sqrt{\overline{e}(t, x)} \mathbf{K}(t, x)}}{n} I_n \]
so for every $(\overline{u}, R)$ such that $\sqrt{\overline{e}(t, x)} \mathbf{K}(t, x) \subset B(\overline{u}, R)$, 
\begin{align*} 
v_0(t, x) \otimes v_0(t, x) - M_0(t, x) &< \frac{|v_0(t, x)|^2 + R^2 - |\overline{u}|^2}{n} I_n \\
&= \frac{|v_0(t, x)|^2 + R^2 - |v_0(t, x) - v_0(t, x) - |\overline{u}|^2}{n} I_n \end{align*}
Therefore, $(v_0(t, x), M_0(t, x), q_0(t, x)) \in ~ \mbox{int} ~ \overline{\mathbf{K}(t, x)(z_0(t, x), \sqrt{\overline{e}(t, x)})}^{co}$ by corollary \ref{cor-enveloppe-convexe}. So $z_0 \in X_0$ which is not empty. 

For $\rho > 0$ and $T > 0$, we set 
\[ J_{T, \rho}(v) = \sup_{t \in [-T, T]} \int_{B_{\rho}} \left[ \overline{e}(t, x) - j_{\mathbf{K}(t, x)}(v(t, x) - v_0(t, x))^2 \right] ~ dx \]
where $B_{\rho} = B(0, \rho)$. 

\begin{Prop} 
\begin{enumerate}
\item If $z = (v, M, q)$ belongs to $X_0$ and $p = q - \frac{|v|^2}{n}$, then $(v, p)$ is solution of the Euler equation with exterior force $f = \nabla_x \cdot \left( F(v) - M \right)$. Furthermore, denoting $(v, M, q) = z_0 + (\widetilde{v}, \widetilde{M}, \widetilde{q})$, we have $\widetilde{v} \in C^{\infty}_c$, $p = q_0 - \frac{|v|^2}{n} + \widetilde{q}$, $\widetilde{q} \in C^{\infty}_c$, and finally $f = \nabla_x \cdot \left( F(v) - M_0 \right) - \nabla_x \widetilde{M}$ and $-\nabla_x \cdot \widetilde{M} =: \widetilde{f} \in C^{\infty}_c$. 
\item Let $(v_k, M_k, q_k)_k$ be a sequence of $X_0$, such that $(v_k, M_k, q_k)_k$ converges in $\mathcal{D}'(\mathbb{R}^{n+1})$ and $(v_k)_k$ converges in $C_{loc}(\mathbb{R}, L^2_{loc}(\mathbb{R}^n))$ to $z = (v, M, q)$, and $J_{T, \rho}(v_k) \to 0$ for every $T, \rho > 0$. Then $v$ is a weak solution of the incompressible Euler equation (with no exterior force), and it satisfies 
\[ \left\{ \begin{array}{l}
v \in C(\mathbb{R}, L^2_w(\mathbb{R}^n, \mathbb{R}^n)) \\
\forall t \in \mathbb{R}, \mbox{ for almost every } x \in \mathbb{R}^n, ~ v(t, x) - v_0(t, x) \in \sqrt{\overline{e}(t, x)} \partial \mathbf{K}(t, x)
\end{array} \right. \]
\end{enumerate} \label{prop-limite-sous-solutions}
\end{Prop}

\begin{Dem}
The first statement is straightforward. 

Then, let $z_k = (v_k, M_k, q_k)$ be as in the second statement, and $(v, M, q)$ its limit. $z_k \in X_0$ implies that $|v_k(t, x) - v_0(t, x)| \leq R_0 \sqrt{\overline{e}(t, x)}$, so $(v_k)$ is uniformly bounded in $L^{\infty}_{loc}(\mathbb{R}, L^2(\mathbb{R}^n, \mathbb{R}^n))$. On the other hand, choosing $\overline{u} = 0$ and $R = R_0 \sqrt{\overline{e}(t, x)}$, we also have 
\begin{align*} M_k(t, x) &\geq v_k(t, x) \otimes v_k(t, x) - \frac{|v_k(t, x)|^2 + R_0^2 \overline{e}(t, x) - |v_k(t, x) - v_0(t, x)|^2}{n} I_n \\
&\geq - \frac{|v_0(t, x)|^2 + 2 |v_k(t, x)-v_0(t, x)| |v_0(t, x)| + R_0^2 \overline{e}(t, x)}{n} I_n \\
&\geq - \frac{(|v_0(t, x)| + R_0 \sqrt{\overline{e}(t, x)})^2}{n} I_n
\end{align*}
Yet $M_k$ being trace-free, we deduce that $|M_k(t, x)| \leq C \left(|v_0(t, x)| + R_0 \sqrt{\overline{e}(t, x)}\right)^2 I_n$ for constant $C > 0$, and so $(M_k)$ is uniformly bounded in $L^{\infty}_{loc}(\mathbb{R}, L^1(\mathbb{R}^n, \mathcal{S}^n_0))$. 

However, 
\[ \partial_t v_k = - \nabla_x \cdot M_k - \nabla_x q_k = \partial_t \mathbb{P} v = - \mathbb{P} \nabla_x \cdot M_k \]
where $\mathbb{P}$ is the Leray projection. This projection acts continuously on $L^2(\mathbb{R}^n)$, so $(\partial_t v_k)$ is uniformly bounded in $L^{\infty}_{loc}(\mathbb{R}, H^{-s}(\mathbb{R}^n, \mathbb{R}^n))$ for $s$ large enough. In particular, by the Aubin-Lions lemma, $(v_k)$ is relatively compact in $C([-T, T], H^{-1}(\mathbb{R}^n, \mathbb{R}^n))$, for every $T > 0$. Up to extracting a subsequence, we therefore obtain convergence in $C(\mathbb{R}, \mathcal{D}') \cap L_{loc}^{\infty}(\mathbb{R}, L^2)$, and so $v \in C(\mathbb{R}, L^2_w)$. 

By linearity, we directly get 
\[ \partial_t v + \nabla_x \cdot M + \nabla_x q = 0, \quad \nabla_x \cdot v = 0 \]

By convexity of $a \mapsto j_{\mathbf{K}(t, x)}(a)^2$ in $\mathbb{R}^n$ ($t, x$ being fixed), we also have $j_{\mathbf{K}(t, x)}(v(t, x) - v_0(t, x))^2 \leq \liminf_k j_{\mathbf{K}(t, x)}(v_k(t, x) - v_0(t, x))^2 \leq \overline{e}(t, x)$, for almost every $x \in \mathbb{R}^n$ and every $t \in \mathbb{R}$. (This inequality comes from Mazur's lemma : we may construct convex combinations of the $v_k - v_0$ strongly converging in $L^2$, so almost everywhere up to extracting a subsequence ; then, at fixed $(t, x)$ such that pointwise convergence holds, we may apply the convexity inequality.) In particular, $v(t, x)-v_0(t, x) \in \sqrt{\overline{e}(t, x)} \mathbf{K}(t, x)$ for almost every $x$ and every $t$. 

Then, by convexity of $(u, N) \mapsto \sigma_{max}(u \otimes u - N)$ the highest eigenvalue, we also deduce that $(v, M, q) \in \overline{\mathbf{K}(t, x)(z_0, \sqrt{\overline{e}(t, x)})}^{co}$ for every $t \in \mathbb{R}$ and almost every $x \in \mathbb{R}^n$. 

Finally, by strong convergence in $C_{loc}(\mathbb{R}, L^2_{loc})$, we have that for every $T, \rho > 0$, 
\[ 0 = \lim J_{T, \rho}(v_k) = J_{T, \rho}(v) \]
hence we deduce that $j_{\mathbf{K}(t, x)}(v(t, x)-v_0(t, x))^2 = \overline{e}(t, x)$ for every $t \in \mathbb{R}$ and almost every $x \in B(0, \rho) \cap \{ \overline{e}(t, \cdot) > 0 \}$. By corollary \ref{cor-enveloppe-convexe}, this implies that $(v, M, q) \in \overline{\mathbf{K}(t, x)(z_0, \overline{e}(t, x))}$ for every $t$ and almost every such $x$. Finally, if $(t, x)$ is a point such that $\overline{e}(t, x) = 0$, we have $j_{\mathbf{K}(t, x)}(v(t, x) - v_0(t, x))^2 = 0 = \overline{e}(t, x)$. In particular, $(v, M, q)(t, x) \in \overline{\mathbf{K}(t, x)(z_0, \sqrt{\overline{e}(t, x)})}$ for every $t$ and almost every $x$, so by proposition \ref{prop-equivalence-systemes} $v$ is a weak solution of the incompressible Euler equation (without exterior force). 
\end{Dem}

\subsection{Improvement step}

\begin{Prop} Let $z = (v, M, q) \in X_0$ and $T, \rho > 0$. There exists a function $\beta : \mathbb{R}^{+*} \to \mathbb{R}^{+*}$ (depending on $T, \rho$) strictly increasing, such that if
\[ J_{T, \rho}(v) \geq \alpha > 0 \]
then for every $\eta > 0$, there exists $z' = (v', M', q') \in X_0$ satisfying
\begin{itemize}
\item $\Vert z - z' \Vert_{C_t H^{-1}_x} \leq \eta$ 
\item $J_{T, \rho}(v') \leq J_{T, \rho}(v) - \beta(\alpha)$
\item $z - z'$ is supported in $[-2T, 2T] \times B_{2\rho}$. 
\end{itemize} \label{prop-amelioration-simple}
\end{Prop}

\begin{Dem}
Without loss of generality, we assume $\alpha \in (0, 1]$. Denote by $\mathcal{L}^n$ the Lebesgue measure of dimension $n$. Moreover, denote by $C_{\rho} > 0$ a constant such that, for every $\varepsilon > 0$ and every regular tiling of $\mathbb{R}^n$ with cubes of side length $\varepsilon$, the number of such cubes intersecting $B_{\rho}$ is bounded by $C_{\rho} \varepsilon^{-n}$. 

For $\varepsilon > 0$ to be chosen small enough later, we define a collection of cubes $C_{\vec i}$ indexed by $\vec i = (i_0, i_1, ..., i_n) \in \mathbb{Z}^{n+1}$ as 
\[ C_{\vec i} = \varepsilon \vec i + \left[ 0, \varepsilon \right)^{n+1} + \frac{1}{2} \mathfrak{d}(\vec i) \varepsilon \vec e_0 \]
where $\mathfrak{d}$ is the parity function defined by : 
\[ \mathfrak{d}(\vec i) = \left\{ \begin{array}{ll}
0 & \mbox{ if } \sum_{j = 1}^n i_j \mbox{ is even,} \\
1 & \mbox{ else.} 
\end{array} \right. \]
and $\vec e_0 = (1, 0, ..., 0)$ is the first basis vector. In particular, in $x$, the $C_{\vec i}$ are a regular tiling of $\mathbb{R}^n$ by cubes of side length $\varepsilon$ ; but the cubes are offset in $t$ every other cube. 

Let then $\mathcal{I}(T, \rho, \varepsilon)$ be the set of $\vec i \in \mathbb{Z}^{n+1}$ such that $C_{\vec i} \cap \left( [-T, T] \times B_{\rho} \right) \neq \emptyset$. For $\varepsilon$ small enough, each of these cubes is included in $[-2T, 2T] \times B_{2\rho}$. 

In every cube $C_{\vec i}$, we define a sub-cube $C_{\vec i}'$ of same center and of side length $\left( \frac{9}{10} \right)^{\frac{1}{n}} \varepsilon$. Denote also by $\overline{e}_{\vec i}$ (respectively $j_{\vec i}$, $\mathbf{K}_{\vec i}, v_{\vec i}$) the value of $\overline{e}$ (respectively $j_{\mathbf{K}}(v-v_0)$, $\mathbf{K}, v - v_0$) at the center of $C_{\vec i}$. By choosing $\varepsilon$ small enough, by uniform continuity of $v, v_0, \overline{e}, \mathbf{K}$ on $[-T, T] \times B_{\rho}$ (compact), we can ensure that the oscillation of $\overline{e} - j_{\mathbf{K}}(v-v_0)^2$ on every cube $C_{\vec i}$ with $\vec i \in \mathcal{I}(T, \rho, \varepsilon)$ is bounded by $\frac{\alpha}{10 \mathcal{L}^n(B_{\rho})}$, and that if for such a $\vec i$ there exists $(t, x) \in C_{\vec i}$ such that $\overline{e}(t, x) < \frac{\alpha}{40 C_{\rho}}$, then $\overline{e} < \frac{\alpha}{30 C_{\rho}}$ on this whole $C_{\vec i}$. 

For every $t_0 \in [-T, T]$, we have moreover that there exists at least one $\epsilon \in \{ 0, 1 \}$ such that for every $\vec i \in \mathbb{Z}^{n+1}$, if $\{ t = t_0 \} \cap C_{\vec i} \neq \emptyset$ and $\mathfrak{d}(\vec i) = \epsilon$, then $\{ t = t_0 \} \cap C_{\vec i}' \neq \emptyset$ : in other words, such a layer in time intersects at least one (equitably distributed) half of these $C_{\vec i}'$. Therefore, by convergence of Riemann sums, we have for every $t_0 \in [-T, T]$ : 
\[ \int_{B_{\rho}} \left[ \overline{e}(t_0, x) - j_{\mathbf{K}(t_0, x)}(v(t_0, x) - v_0(t_0, x))^2 \right] dx \leq 3 \sum_{\substack{\vec i \in \mathcal{I}(T, \rho, \varepsilon), \\ C'_{\vec i} \cap \{ t = t_0 \} \neq \emptyset}} \mathcal{L}^n(C'_{\vec i} \cap \{ t = t_0 \}) ~ \left[ \overline{e}_{\vec i} - j_{\vec i}^2 \right] \quad + \frac{\alpha}{10} \]

Let $c > 0$ to be fixed later. The contribution in the above right-hand sum of cubes such that $\overline{e}_{\vec i} - j_{\vec i}^2 \leq c \alpha$ is bounded by 
\[ C_{\rho} \varepsilon^{-n} * \frac{9}{10} \varepsilon^n c \alpha \leq C_{\rho} c \alpha \]
In particular, if $c = \frac{1}{30 C_{\rho}}$, we get
\[ \int_{B_{\rho}} \left[ \overline{e}(t_0, x) j_{\mathbf{K}(t_0, x)}(v(t_0, x) - v_0(t_0, x))^2 \right] ~ dx \leq 3 \sum_{\substack{\vec i \in \mathcal{I}(T, \rho, \varepsilon), \\ C'_{\vec i} \cap \{ t = t_0 \} \neq \emptyset, \\ \overline{e}_{\vec i} - j_{\vec i}^2 > c \alpha}} \mathcal{L}^n(C'_{\vec i} \cap \{ t = t_0 \}) ~ \left[ \overline{e}_{\vec i} - j_{\vec i}^2 \right] \quad + \frac{\alpha}{5} \]
Moreover, for any $\vec i$ in the sum, $\overline{e}(t, x) > 0$ for every $(t, x) \in C_{\vec i}$, and so $\sqrt{\overline{e}(t, x)} \mathbf{K}(t, x)$ is never reduced to $\{ 0 \}$. 

In particular, if $t_0$ is such that 
\[ \int_{B_{\rho}} \left[ \overline{e}(t_0, x) - j_{\mathbf{K}(t_0, x)}(v(t_0, x)-v_0(t_0, x))^2 \right] ~ dx \geq \frac{4 \alpha}{5} \]
we have that
\[ \sum_{\substack{\vec i \in \mathcal{I}(T, \rho, \varepsilon), \\ C_{\vec i}' \cap \{ t = t_0 \} \neq \emptyset, \\ \overline{e}_{\vec i} - j_{\vec i}^2 > c \alpha}} \mathcal{L}^n(C_{\vec i}' \cap \{ t = t_0 \}) ~ \left[ \overline{e}_{\vec i} - j_{\vec i}^2 \right] \quad \geq \frac{\alpha}{5} \]

We will now construct $(v', M', q')$ as $(v', M', q') = (v, M, q) + \sum_{\vec i} (\widetilde{v_{\vec i}}, \widetilde{M_{\vec i}}, \widetilde{q_{\vec i}})$, satisfying the following properties : 
\begin{itemize}
\item every $\widetilde{z_{\vec i}} = (\widetilde{v_{\vec i}}, \widetilde{M_{\vec i}}, \widetilde{q_{\vec i}})$ is a solution of \eqref{equ-Euler-lin} and is supported in $C_{\vec i}$
\item if $\vec i$ is not in the above sum, $\widetilde{z_{\vec i}} = 0$
\end{itemize}

In particular, the gain in the functional $J_{T, \rho}$ at any time $t_0$ will be : 
\[ \begin{aligned} 
\Delta(t_0) &= \int_{B_{\rho}} \left[ \overline{e}(t_0, x) - j_{\mathbf{K}(t_0, x)}(v(t_0, x) - v_0(t_0, x))^2 \right] ~ dx \\
&\quad - \int_{B_{\rho}} \left[ \overline{e}(t_0, x) - j_{\mathbf{K}(t_0, x)}(v'(t_0, x) - v_0(t_0, x))^2 \right] ~ dx \\
&= \int_{B_{\rho}} \left[ j_{\mathbf{K}(t_0, x)}(v'(t_0, x) - v_0(t_0, x))^2 - j_{\mathbf{K}(t_0, x)}(v(t_0, x) - v_0(t_0, x))^2 \right] ~ dx \\
&= \sum_{\substack{ \vec i \in \mathcal{I}(T, \rho, \varepsilon), \\ C_{\vec i}' \cap \{ t = t_0 \} \neq \emptyset, \\ \overline{e}_{\vec i} - j_{\vec i}^2 > c \alpha }} \int_{C_{\vec i} \cap \{ t = t_0 \}} \Big[ j_{\mathbf{K}(t_0, x)}(v(t_0, x) + \widetilde{v_{\vec i}}(t_0, x) - v_0(t_0, x))^2 \\
&\quad \quad \quad \quad \quad \quad \quad \quad \quad \quad \quad \quad \quad \quad \quad \quad \quad \quad - j_{\mathbf{K}(t_0, x)}(v(t_0, x) - v_0(t_0, x))^2 \Big] ~ dx 
\end{aligned} \]
For each $\vec i$, constructing the subsolution $\widetilde{z_{\vec i}}$ by proposition \ref{prop-solutions-loc} (or setting $\widetilde{z_{\vec i}} = 0$) with a $\lambda_{\vec i} = (\overline{v}_{\vec i}, \overline{M}_{\vec i}, \overline{q}_{\vec i})$ to be determined, one has 
\[ \Delta(t_0) \geq \sum_{\substack{ \vec i \in \mathcal{I}(T, \rho, \varepsilon), \\ C_{\vec i}' \cap \{ t = t_0 \} \neq \emptyset, \\ \overline{e}_{\vec i} > c \alpha }} \frac{c_0 r_0^2}{6} \mathcal{L}^n(C_{\vec i}' \cap \{ t = t_0 \}) |\overline{v}_{\vec i}|^2 \]
Note that this is always non-negative. 

We may now distinguish two cases for $t_0$. If $t_0$ satisfies  
\[ \int_{B_{\rho}} \left[ \overline{e}(t_0, x)  - j_{\mathbf{K}(t_0, x)}(v(t_0, x)-v_0(t_0, x))^2 \right] dx < \frac{4 \alpha}{5} \]
then this bound remains replacing $v$ by $v'$, since $\Delta(t_0) \geq 0$. 

Conversely, if $t_0$ is such that 
\[ \int_{B_{\rho}} \left[ \overline{e}(t_0, x) - j_{\mathbf{K}(t_0, x)}(v(t_0, x)-v_0(t_0, x))^2 \right] dx \geq \frac{4 \alpha}{5} \]
then we saw that 
\[ \sum_{\substack{\vec i \in \mathcal{I}(T, \rho, \varepsilon), \\ 
C_{\vec i}' \cap \{ t = t_0 \} \neq \emptyset, \\
\overline{e}_{\vec i} - j_{\vec i}^2 > c \alpha}} \mathcal{L}^n(C'_{\vec i} \cap \{ t = t_0 \}) \left[ \overline{e}_{\vec i} - j_{\vec i}^2 \right] \quad \geq \frac{\alpha}{5} \]
In particular, if we manage to get 
\[ |\overline{v}_{\vec i}(t_0, x)|^2 \geq \frac{\beta}{\alpha} \left( \overline{e}_{\vec i} - j_{\vec i}^2 \right) \]
for some $\beta = \beta(\alpha)$, then we will obtain 
\[ \begin{aligned}
\Delta(t_0) &\geq \sum_{\substack{\vec i \in \mathcal{I}(T, \rho, \varepsilon), \\
C_{\vec i}' \cap \{ t = t_0 \} \neq \emptyset, \\
\overline{e}_{\vec i} - j_{\vec i}^2 > c \alpha}} \frac{c_0 r_0^2}{6} \mathcal{L}^n(C'_{\vec i} \cap \{ t = t_0 \}) \frac{\beta}{\alpha} \left[ \overline{e}_{\vec i} - j_{\vec i}^2 \right] \\
&\geq \frac{c_0 r_0^2 \beta}{30}
\end{aligned} \]
and $\frac{c_0 r_0^2}{30}$ is a constant independant of $\alpha, \varepsilon$. 

In both cases, we manage to strictly decrease the functional $J_{T, \rho}$. It is now enough to construct such a $\widetilde{z_{\vec i}}$ on every involved cube. 

Let
\[ \varrho = \inf \left\{ \mbox{dist}\left( z(x, t), \partial \left( \overline{\mathbf{K}(t, x, z_0(t, x), \sqrt{\overline{e}(t, x)})}^{co} \right) \right), ~ (t, x) \in [-T, T] \times B_{\rho}, ~ \overline{e}(t, x) \geq \frac{\alpha}{40 C_{\rho}} \right\} > 0 \]
defined independently of $\varepsilon$. We can then construct by corollary \ref{cor-enveloppe-convexe}, for every $\vec i$ such that $\overline{e}_{\vec i} - j_{\vec i}^2 > c \alpha$, a $\lambda_{\vec i} = (\overline{v}_{\vec i}, \overline{M}_{\vec i}, \overline{q}_{\vec i})$, direction of oscillation such that $z_{\vec i} + [-\lambda_{\vec i}, \lambda_{\vec i}]$ is included in the interior of $\overline{\mathbf{K}_{\vec i}(z_{0, \vec i}, \sqrt{\overline{e}_{\vec i}})}^{co}$, and at a distance at least $\frac{\varrho}{2}$ of the boundary, and 
\[ |\overline{v}_{\vec i}| \geq C(n) \mbox{dist}(v_{\vec i}, \sqrt{\overline{e}_{\vec i}} \partial \mathbf{K}_{\vec i}) \]
Let then $w_{\vec i} \in \sqrt{\overline{e}_{\vec i}} \partial \mathbf{K}_{\vec i}$ be the point realizing this distance : by remark \ref{rem-equiv-jauge}, 
\begin{align*} 
|\overline{v}_{\vec i}| &\geq C(n) |v_{\vec i} - w_{\vec i}| \\
&\geq C(n) r_0 j_{\mathbf{K}_{\vec i}}(w_{\vec i}-v_{\vec i}) \\
&\geq C(n) r_0 \left( j_{\mathbf{K}_{\vec i}}(w_{\vec i}) - j_{\mathbf{K}_{\vec i}}(v_{\vec i}) \right) \\
&\geq \frac{C(n) r_0}{2 \sqrt{\overline{e}_{\vec i}}} \left( \overline{e}_{\vec i} - j_{\mathbf{K}_{\vec i}}(v_{\vec i})^2 \right) 
\end{align*}
Now $\overline{e}$ is bounded on $[-T, T] \times B_{\rho}$ by continuity, independently of $\alpha, \varepsilon$, and we only consider $\vec i$ such that $\overline{e}_{\vec i} - j_{\vec i}^2 \geq c \alpha$, so there exists a constant $C = C(n, T, \rho, \overline{e}, r_0)$ such that
\[ |\overline{v}_{\vec i}^2 \geq C(n, T, \rho, \overline{e}, r_0) \alpha \left( \overline{e}_i - j_{\vec i}^2 \right) \]
So we may choose $\beta(\alpha) = C(n, T, \rho, \overline{e}, r_0) \alpha^2$. 

Up to choosing $\varepsilon$ small enough, by uniform continuity of $z$ and of $\mathbf{K}$, we may keep $z(t, x) + [-\lambda_{\vec i}, \lambda_{\vec i}]$ in the interior of $\overline{\mathbf{K}(t, x, z_0(t, x), \sqrt{\overline{e}(t, x)})}^{co}$ for every $(t, x) \in C_{\vec i}$. We thus find a neighborhood $\mathcal{V}_{\vec i}$ of $[-\lambda_{\vec i}, \lambda_{\vec i}]$ such that $z(t, x) + \mathcal{V}_{\vec i} \subset ~ \mbox{int} ~ \overline{\mathbf{K}(t, x, z_0(t, x), \sqrt{\overline{e}(t, x)})}^{co}$ for every $(t, x) \in C_{\vec i}$. When we apply proposition \ref{prop-solutions-loc} to construct a localized wave $\widetilde{z}_{\vec i}$, we can take it supported in $C_{\vec i}$, valued in $\mathcal{V}$, and we may choose the $C_tH^{-1}_x$ norm of every added subsolution to be arbitrarily small. Since at this point we have already fixed the size of the paving, and hence the maximal number of such added subsolutions, we may satisfy as well the point 3 of proposition \ref{prop-amelioration-simple}, and this concludes the proof. 
\end{Dem}

\begin{Prop} Let $z = (v, M, q) \in X_0$, $(T_l)$ and $(\rho_l)$ two strictly increasing sequences, $l = 1, ..., L$. Assume that 
\[ J_{T_l, \rho_l}(v) \geq \alpha_l > 0 \]
and let $\eta > 0$. Then there exists $z' = (v', M', q') \in X_0$ such that 
\begin{itemize}
\item $\Vert z - z' \Vert_{C_t H^{-1}_x} \leq \eta$ 
\item $\forall ~ l \in \{ 1, ..., L \}$, $J_{T_l, \rho_l}(v') \leq J_{T_l, \rho_l}(v) - \beta_l(\alpha_l)$
\end{itemize}
where $\beta_l = \beta_{T_l, \rho_l}$ is the function of proposition \ref{prop-amelioration-simple}. \label{prop-amelioration-L}
\end{Prop}

\begin{Dem}
We proceed as in the proof of proposition \ref{prop-amelioration-simple}. The beginning remains the same, and we get all the bounds for every $l$ (uniformly), and in particular, if $t_0 \in \left[ -T_l, T_l \right]$ is such that
\[ \int_{B_{\rho_l}} \left[ \overline{e}(t_0, x) - j_{\mathbf{K}(t_0, x)}(v(t_0, x)-v_0(t_0, x))^2 \right] dx \geq \frac{3 \alpha_l}{5}, \] 
then we have
\[ \sum_{\substack{\vec i \in \mathcal{I}(T_l, \rho_l, \varepsilon),\\ 
C_{\vec i} \cap \{ t = t_0 \} \neq \emptyset, \\
\overline{e}_{\vec i} > c \alpha_l}} \mathcal{L}^n(C'_{\vec i} \cap \{ t = t_0 \}) \overline{e}_{\vec i} \geq \frac{\alpha_l}{5} \]

Then, we may not necessarily keep the condition $\widetilde{z_{\vec i}} = 0$ when $\overline{e}_{\vec i} \leq c \alpha_l$, but we may refine the computation of $\Delta_l(t_0)$ : 
\begin{align*}
\Delta_l(t_0) &= \int_{B_{\rho_l}} \left[ \overline{e}(t_0, x) - j_{\mathbf{K}(t_0, x)}(v(t_0, x)-v_0(t_0, x))^2 \right] dx \\
&\quad - \int_{B_{\rho_l}} \left[ \overline{e}(t_0, x) - j_{\mathbf{K}(t_0, x)}(v'(t_0, x)-v_0(t_0, x))^2 \right] dx \\
&= \int_{B_{\rho_l}} \left( j_{\mathbf{K}(t_0, x)}(v'(t_0, x)-v_0(t_0, x))^2 - j_{\mathbf{K}(t_0, x)}(v(t_0, x)-v_0(t_0, x))^2 \right) dx \\
&= \sum_{\substack{\vec i, \\
C_{\vec i}' \cap \left( \{ t_0 \} \times B_{\rho_l} \right) \neq \emptyset, \\
\overline{e}_{\vec i} - j_{\vec i}^2 > c \min_l \alpha_l}} \int_{C_{\vec i} \cap \{ t_0 \} \times B_{\rho_l}} \Bigl[ j_{\mathbf{K}(t_0, x)}(v(t_0, x) + \widetilde{v_{\vec i}}(t_0, x)-v_0(t_0, x))^2 \\
&\quad \quad \quad \quad \quad \quad \quad \quad \quad \quad \quad \quad \quad \quad \quad \quad \quad \quad - j_{\mathbf{K}(t_0, x)}(v(t_0, x)-v_0(t_0, x))^2 \Bigl] dx \\
&\\
&\geq \sum_{\substack{\vec i, \\
C_{\vec i}' \cap \left( \{ t_0 \} \times B_{\rho_l} \right) \neq \emptyset, \\
\overline{e}_{\vec i} - j_{\vec i}^2 > c \min_l \alpha_l}} \frac{c_0 r_0^2}{6} \mathcal{L}^n(C_{\vec i}' \cap \{ t = t_0 \}) |\overline{v}_{\vec i}|^2 \\
&\geq \sum_{\substack{\vec i \in \mathcal{I}(T_l, \rho_l, \varepsilon) \\
C_{\vec i}' \cap \{ t = t_0 \} \neq \emptyset, \\
\overline{e}_{\vec i} - j_{\vec i}^2 \geq c \alpha_l}} \frac{c_0 r_0^2}{6} \mathcal{L}^n(C_{\vec i}' \cap \{ t = t_0 \}) |\overline{v}_{\vec i}|^2 
\end{align*}
as before. The end of the proof is the same. 
\end{Dem}

\subsection{Iteration}

\begin{Dem}[of theorem \ref{theo-principal}]
Let $\psi \in C^{\infty}_c(\mathbb{R}^n, \mathbb{R}^{+})$ and $\chi \in C^{\infty}_c(\mathbb{R}, \mathbb{R}^{+})$ be of total mass $1$, supported in the unit ball and in $[-1, 1]$ respectively. Set  
\begin{equation} \rho_{\varepsilon}(t, x) = \varepsilon^{-n-1} \psi\left( \frac{x}{\varepsilon} \right) \chi\left( \frac{t}{\varepsilon} \right) = \psi_{\varepsilon}(x) \chi_{\varepsilon}(t) \label{equ-def-rho-convolution} \end{equation}
a regularizing kernel in both space and time. We denote by $\ast$ the convolution product in space-time on $\mathbb{R}^{n+1}$, $\ast_x$ the convolution product in $x$ only, and $\ast_t$ in $t$ only. 

Denote by $\beta_j = \beta_{j, j}$ the function of propositions \ref{prop-amelioration-simple} and \ref{prop-amelioration-L} for $T = j, \rho = j$, with $j > 0$. 

We start from $z_0$ which lies in $X_0$, and we set $\varepsilon_0 = 1$. 

Assume that $(z_j, \varepsilon_j) \in X_0 \times \mathbb{R}^{+*}$ is constructed for $j \leq k-1$, and apply proposition \ref{prop-amelioration-L} to construct $z_k = (v_k, M_k, q_k)$ such that 
\begin{itemize}
\item $J_{j, j}(v_k) \leq J_{j, j}(v_{k-1}) - \beta_j\left( J_{j, j}(v_{k-1}) \right)$ for every $j \leq k$
\item $\sup_t \Vert z_k(t) - z_{k-1}(t) \Vert_{H^{-1}} \leq \eta 2^{-k}$ 
\item and for every $j \leq k-1$, $\sup_t \Vert (z_k - z_{k-1}) \ast \rho_{\varepsilon_j} \Vert_{L^2} \leq 2^{-k}$. 
\end{itemize}
The two first points follow directly from proposition \ref{prop-amelioration-L}. For the last one, we have 
\begin{align*}
\sup_t \Vert \left( (z_k - z_{k-1}) \ast \rho_{\varepsilon_j} \right)(t) \Vert_{L^2}
&= \sup_t \left\Vert \left( \left( (z_k - z_{k-1}) \ast_t \chi_{\varepsilon_j} \right) \ast_x \psi_{\varepsilon_j} \right)(t) \right\Vert_{L^2} \\
&= \sup_t \left\Vert \mathcal{F} \left( (z_k - z_{k-1}) \ast_t \chi_{\varepsilon_j} \right) (t) \mathcal{F}(\psi_{\varepsilon_j}) \right\Vert_{L^2} \\
& \quad \quad \quad \quad \quad \quad \quad \quad \quad \quad \quad \quad \mbox{ (Parseval's identity)} \\
&\leq \sup_t \left\Vert \langle \xi \rangle^{-1} \left( \mathcal{F}(z_k - z_{k-1}) \ast_t \chi_{\varepsilon_j} \right)(t) \right\Vert_{L^2(\mathbb{R}^n)} \Vert \langle \xi \rangle \mathcal{F} (\psi_{\varepsilon_j}) \Vert_{L^{\infty}} \\
& \quad \quad \quad \quad \quad \quad \quad \quad \quad \quad \quad \quad \mbox{ (Hölder's inequality)} \\
&\leq C \varepsilon_j^{-1} \sup_t \left( \left\Vert \langle \xi \rangle^{-1} \mathcal{F}(z_k - z_{k-1}) \right\Vert_{L^2} \ast_t \chi_{\varepsilon_j} \right)(t) \\
& \quad \quad \quad \quad \quad \quad \quad \quad \quad \quad \quad \quad  \mbox{ (integral Minkowski's inequality)} \\
&\leq C \varepsilon_j^{-1} \sup_t \Vert z_k - z_{k-1} \Vert_{H^{-1}} \\
& \quad \quad \quad \quad \quad \quad \quad \quad \quad \quad \quad \quad  \mbox{ (Parseval's identity)}
\end{align*}
In particular, we may satisfy this last point up to choosing $\Vert z_k - z_{k-1} \Vert_{H^{-1}}$ small enough. 

Finally, we choose $\varepsilon_k \in (0, 1)$ such that
\[ \forall j \leq k, \quad \sup_t \Vert z_j - z_j \ast \rho_{\varepsilon_k} \Vert_{L^2} \leq 2^{-k} \]

By construction, the sequence $(z_k)$ is a Cauchy sequence in $L^{\infty}_{t, loc} H^{-1}_x$, hence it converges weakly in the sense of distributions to a $z = (v, M, q) \in L^{\infty}_{t, loc} H^{-1}_x$, with $\Vert z - z_0 \Vert_{L^{\infty}_t H^{-1}_x} \leq \eta$ (globally in time). Furthermore, since $z_k \in X_0$ for every $k$, $v_k \in L^{\infty}_{t, loc} L^2_x$ uniformly, we get $v \in L^{\infty}_{t, loc} L^2_x$. $\rho_{\varepsilon_k}$ being a regularization kernel, we have convergence of $z_l \ast \rho_{\varepsilon_k}$ to $z \ast \rho_{\varepsilon_k}$ as $l \to \infty$ ($k$ being fixed) in $C_t L^2_x$ for instance. 

Let now $j \geq 1$ and $E = C([-j, j], L^2(B_j))$. For $k \geq j$, 
\[ \Vert z_k - z \Vert_E \leq \Vert z_k - z_k \ast \rho_{\varepsilon_k} \Vert_E + \Vert z_k \ast \rho_{\varepsilon_k} - z \ast \rho_{\varepsilon_k} \Vert_E + \Vert z \ast \rho_{\varepsilon_k} - z \Vert_E \]
The last term converges to $0$ as $k \to \infty$. The choice of $\varepsilon_k$ also ensures that the first term is bounded by $2^{-k}$, so converges to $0$. Finally, for the second term : 
\[ \Vert z_k \ast \rho_{\varepsilon_k} - z \ast \rho_{\varepsilon_k} \Vert_E = \lim_{l \to \infty} \Vert z_k \ast \rho_{\varepsilon_k} - z_l \ast \rho_{\varepsilon_k} \Vert_E \leq \sum_{m \geq k} \Vert (z_m - z_{m+1}) \ast \rho_{\varepsilon_k} \Vert_E \leq 2^{-k} \]
by construction. This implies local uniform convergence of the velocity. 

Let now $T, \rho > 0$. There exists $j$ such that $\max(T, \rho) \leq j$. In particular, $J_{T, \rho} \leq J_{j, j}$. 

Yet $J_{l, l}(v_k)$ converges to $0$ for every (fixed) $l$. Indeed, if $\alpha_l = \limsup_{k \to \infty} J_{l, l}(v_k) > 0$, then by passing to the limit in the improvement inequality we have $\alpha_l \leq \alpha_l - \beta_l(\alpha_l) < \alpha_l$, which is a contradiction. Therefore this upper limit is $0$, hence the convergence. 

We conclude that all the hypothesis of the second statement of proposition \ref{prop-limite-sous-solutions} are satisfied, which ends the proof of theorem \ref{theo-principal}. 
\end{Dem}

\section{Construction of continuous paths} \label{section-chemins}

\begin{Lem} Let $u, v \in C^0_t L^2_x$, and $\rho_{\varepsilon}$ the approximation of unity as in \eqref{equ-def-rho-convolution}. Then for every $T > 0$, 
\[ \Vert (uv) \ast \rho_{\varepsilon} - (u \ast \rho_{\varepsilon}) (v \ast \rho_{\varepsilon}) \Vert_{L^{\infty}([-T, T], L^2_x)} \underset{\varepsilon \to 0}{\longrightarrow} 0 \] \label{lem-convergence-approx-cunif}
\end{Lem}

\begin{Dem}
It is enough to show convergence in $L^{\infty}_t L^2_x$ of $(uv) \ast \rho_{\varepsilon}$ to $uv$, and of $(u \ast \rho_{\varepsilon}) (v \ast \rho_{\varepsilon})$ to $uv$. 

Let $\delta > 0$. There exists $\varepsilon > 0$ such that, for every $s, t \in [-T-1, T+1]$ satisfying $|t-s| < \varepsilon$ we have
\[ \Vert u(t) - u(s) \Vert_{L^2} + \Vert v(t) - v(s) \Vert_{L^2} < \delta \]
by uniform continuity on $[-T-1, T+1]$. In particular, for every $t$, 
\[ \Vert u(t) - u \ast \rho_{\varepsilon}(t) \Vert_{L^2} \leq \sup_{|s - t| < \varepsilon} \Vert u(t) - u(s) \Vert_{L^2} \leq \delta \]
and likewise for $v$. We deduce that $u \ast \rho_{\varepsilon}$ converges to $u$ in $L^{\infty}_t L^2_x$, and likewise for $v$, then that $(u \ast \rho_{\varepsilon}) (v \ast \rho_{\varepsilon})$ converges to $uv$ in $L^{\infty}_t L^1_x$. 

By noticing that $uv$ is in $C^0_t L^1_x$ and uniformly continuous on $[-T-1, T+1]$, and that the above argument for $L^2$ is the same for any $L^p$, we can show the same convergence for $(uv) \ast \rho_{\varepsilon}$. 
\end{Dem}

We now prove theorem \ref{theo-cheminsCdemi}. To construct a path, we start by contructing intermediate solutions between two given solutions $u_0$ and $u_1$. 

\begin{Lem} Let $(u_0, p_0)$ and $(u_1, p_1)$ be two weak solutions of the incompressible Euler equation, with $u_i \in C^0_t L^2_x$, $i = 1, 2$. Let $\delta > 0, T > 0, T' > 0$. There exists $(u_{1/2}, p_{1/2})$ a weak solution of the incompressible Euler equation such that 
\[ \forall ~ 0 < t \leq T, \quad \sup_{i = 1, 2} \Vert u_{1/2} - u_i \Vert_{L^{\infty}([-t, t], L^2_x)} \leq \frac{\Vert u_1 - u_0 \Vert_{L^{\infty}([-t-T', t+T'], L^2_x)}}{\sqrt{2}} + \delta \]
et $u_{1/2} \in C^0_t L^2_x$. \label{lem-construction-intermediaire} 
\end{Lem}

\begin{Dem}
Let $(u_0, p_0)$ and $(u_1, p_1)$ be as in the lemma. Define 
\[ z_i = \left( u_i, F(u_i), p_i + \frac{|u_i|^2}{n} \right) \]

Let $\varepsilon > 0$, and set
\[ z_i^{\varepsilon} := \rho_{\varepsilon} \ast z_i \]
Denote $z_i^{\varepsilon} = (v_i^{\varepsilon}, M_i^{\varepsilon}, q_i^{\varepsilon})$. Here above, $\rho_{\varepsilon}$ is the same regularization kernel as in the proof of theorem \ref{theo-principal}, see formula \eqref{equ-def-rho-convolution} (so it regularizes in both space and time). Then $v_i^{\varepsilon} \in L^{\infty}_t L^2_x$, and for every $\varepsilon > 0$, $z_i^{\varepsilon} \in C^0_{t, x}$, with a bound that degenerates as $\varepsilon \to 0$. Moreover, $z_i^{\varepsilon}$ is solution of \eqref{equ-Euler-lin} for every $\varepsilon$ by linearity of this equation. We may then define 
\[ z_{1/2}^{\varepsilon} = \frac{1}{2} (z_0^{\varepsilon} + z_1^{\varepsilon}) = \left( u_{1/2}^{\varepsilon}, M_{1/2}^{\varepsilon}, q_{1/2}^{\varepsilon} \right) \]
that satisfies the same properties. 

We now apply theorem \ref{theo-principal} from $z_{1/2}^{\varepsilon}$ to construct $u_{1/2}$ a weak solution of the incompressible Euler equation for which we will control the distance to $u_0$ and $u_1$ simultaneously, by choosing suitable convex compact sets $\mathbf{K}(t, x)$. 

For this, we define
\[ \mathbf{K}^{0, \varepsilon}(t, x) = B\left( u_0^{\varepsilon}(t, x), R^{\varepsilon}(t, x) \right) \cap B\left( u_1^{\varepsilon}(t, x), R^{\varepsilon}(t, x) \right) - u_{1/2}^{\varepsilon}(t, x) \]
with a radius $R^{\varepsilon}(t, x)$ to be chosen suitably in $C^0_{t, x} \cap L^{\infty}_{t, loc} L^2_x$, and $R^{\varepsilon}(t, x) > \frac{\sqrt{2}}{2} |u_0^{\varepsilon}(t, x) - u_1^{\varepsilon}(t, x)|$. It is then clear that $\mathbf{K}^{0, \varepsilon}(t, x) \in \mathfrak{K}$ for every $t, x$. Furthermore, the map $\mathbf{K}^{0, \varepsilon}$ is continuous for the Hausdorff distance. We may define the renormalized convex as
\[ \mathbf{K}^{\varepsilon}(t, x) = \frac{1}{r_{max}(\mathbf{K}^{0, \varepsilon}(t, x))} \mathbf{K}^{0, \varepsilon}(t, x) \]
and set
\[ \overline{e}^{\varepsilon}(t, x) = r_{max}(\mathbf{K}^{0, \varepsilon}(t, x)) \] 
$\overline{e}^{\varepsilon}$ is automatically continuous by lemma \ref{lem-continuite-rminmax}, and it is bounded by $R^{\varepsilon}(t, x)^2$, so it also belongs to $L^{\infty}_{t, loc} L^1_x$. Furthermore, that way, 
\[ r_{max}(\mathbf{K}^{\varepsilon}(t, x)) = 1 \]
so we have $\mathbf{K}^{\varepsilon}(t, x) \subset B(0, 1)$ for any $t, x$. 

By remark \ref{rem-unifcvxte} and lemma \ref{lem-uniforme-convexite-global}, we also know that $\mathbf{K}^{\varepsilon}$ is $c_0$-uniformly convex for a constant $c_0 > 0$, and that $r_{min}(\mathbf{K}^{\varepsilon}(t, x)) \geq r_0$ for some $r_0 > 0$, which means $B(0, r_0) \subset \mathbf{K}^{\varepsilon}(t, x)$ for every $t, x$. 

By lemma \ref{lem-calcul-aK}, we have
\[ \overline{e}^{\varepsilon}(t, x) a_{\mathbf{K}^{\varepsilon}(t, x)} = R^{\varepsilon}(t, x)^2 - \frac{1}{4} |u_1^{\varepsilon}(t, x) - u_0^{\varepsilon}(t, x)|^2 \]
So to satisfy the last hypothesis of theorem \ref{theo-principal}, we should choose $R(t, x)$ such that
\[ R^{\varepsilon}(t,x)^2 I_n > \frac{1}{4} |u_1^{\varepsilon}(t, x) - u_0^{\varepsilon}(t, x)|^2 I_n + n \left( F(v_{1/2}^{\varepsilon}(t, x)) - M_{1/2}^{\varepsilon}(t, x) \right) \]
in the sense of matrices. 

Let us first compute this for $\varepsilon = 0$. Then, we need : 
\[ \begin{aligned}
R(t, x)^2 I_n &> \frac{1}{4} |u_1(t, x) - u_0(t, x)|^2 I_n \\
&\quad + n \Biggl( \left( \frac{u_0(t, x) + u_1(t, x)}{2} \right) \otimes \left( \frac{u_0(t, x) + u_1(t, x)}{2} \right) - \frac{|u_0(t, x) + u_1(t, x)|^2}{4n} I_n \\
&\quad \quad \quad - \frac{u_0(t, x) \otimes u_0(t, x) + u_1(t, x) \otimes u_1(t, x)}{2} + \frac{|u_0(t, x)|^2 + |u_1(t, x)|^2}{2n} I_n \Biggl) \\
&= \frac{|u_1(t, x) - u_0(t, x)|^2}{4} I_n \\
&\quad + n \left( - \left( \frac{u_1(t, x) - u_0(t, x)}{2} \right) \otimes \left( \frac{u_1(t, x) - u_0(t, x)}{2} \right) + \frac{|u_1(t, x) - u_0(t, x)|^2}{4n} I_n \right) \\
&= \frac{|u_1(t, x) - u_0(t, x)|^2}{2} I_n - n \left( \frac{u_1(t, x) - u_0(t, x)}{2} \right) \otimes \left( \frac{u_1(t, x) - u_0(t, x)}{2} \right)
\end{aligned} \]
which is equivalent to 
\[ R(t, x)^2 > \frac{|u_1(t, x) - u_0(t, x)|^2}{2} \]
which is exactly the hypothesis $R(t, x) > \frac{\sqrt{2}}{2} |u_1(t, x) - u_0(t, x)|$. 

In general, we want 
\[ \begin{aligned}
R^{\varepsilon}(t, x)^2 I_n &> \frac{1}{4} |u_1^{\varepsilon}(t, x) - u_0^{\varepsilon}(t, x)|^2 I_n + n \left( F(v_{1/2}^{\varepsilon}(t, x)) - M_{1/2}^{\varepsilon}(t, x) \right) \\
&= \frac{|u_1^{\varepsilon}(t, x) - u_0^{\varepsilon}(t, x)|^2}{2} I_n - n \left( \frac{u_1^{\varepsilon}(t, x) - u_0^{\varepsilon}(t, x)}{2} \right) \otimes \left( \frac{u_1^{\varepsilon}(t, x) - u_0^{\varepsilon}(t, x)}{2} \right) \\
&\quad \quad \quad + \frac{\left( |u_0|^2 \right)^{\varepsilon}(t, x) - |u_0^{\varepsilon}(t, x)|^2 + \left( |u_1|^2 \right)^{\varepsilon}(t, x) - |u_1^{\varepsilon}(t, x)|^2}{2} I_n \\
&\quad \quad \quad - n \frac{(u_0 \otimes u_0)^{\varepsilon}(t, x) - u_0^{\varepsilon}(t, x) \otimes u_0^{\varepsilon}(t, x) + (u_1 \otimes u_1)^{\varepsilon}(t, x) - u_1^{\varepsilon}(t, x) \otimes u_1^{\varepsilon}(t, x)}{2}
\end{aligned} \]
where we extend the notation $( \cdot )^{\varepsilon} := ( \cdot ) \ast \rho_{\varepsilon}$, using the same computation as in the case $\varepsilon = 0$. Therefore, we can choose
\[ R^{\varepsilon}(t, x)^2 = \frac{|u_1^{\varepsilon}(t, x) - u_0^{\varepsilon}(t, x)|^2}{2} + \tau_{\varepsilon}(t, x) \]
for a remainder $\tau_{\varepsilon}$ converging to $0$ in $L^{\infty}_{loc, t} L^2_x$ as $\varepsilon \to 0$ by lemma \ref{lem-convergence-approx-cunif}. 

We can now apply theorem \ref{theo-principal}. This gives us $(u_{1/2}, p_{1/2})$ a solution of the incompressible Euler equation, with
\[ \left| u_{1/2}(t, x) - u_i^{\varepsilon}(t, x) \right| \leq R^{\varepsilon}(t, x) \]
for $i = 1, 2$. In particular, for every $T$, 
\[ \Vert u_{1/2} - u_i^{\varepsilon} \Vert_{L^{\infty}([-T, T], L^2_x)} \leq \Vert R^{\varepsilon} \Vert_{L^{\infty}([-T, T], L^2_x)} \leq \frac{\Vert u_1^{\varepsilon} - u_0^{\varepsilon} \Vert_{L^{\infty}([-T, T] L^2_x)}}{\sqrt{2}} + \Vert \tau_{\varepsilon} \Vert_{L^{\infty}([-T, T], L^1_x)} \]

We then notice that 
\[ \Vert u_1^{\varepsilon} - u_0^{\varepsilon} \Vert_{L^{\infty}([-T, T], L^2_x)} = \Vert (u_1 - u_0) \ast \rho_{\varepsilon} \Vert_{L^{\infty}([-T, T], L^2_x)} \leq \Vert u_1 - u_0 \Vert_{L^{\infty}([-T-T', T+T'], L^2_x)} \]
by Young's convolution inequality, up to choosing $\varepsilon < T'$, and we saw that $\Vert \tau_{\varepsilon} \Vert_{L^{\infty}([-T, T], L^1_x)}$ can be made arbitrarily small by choosing $\varepsilon$ small enough. 

Finally, 
\[ \Vert u_{1/2} - u_i \Vert_{L^{\infty}([-T, T], L^2_x)} \leq \Vert u_{1/2} - u_i^{\varepsilon} \Vert_{L^{\infty}([-T, T], L^2_x)} + \Vert u_i^{\varepsilon} - u_i \Vert_{L^{\infty}([-T, T], L^2_x)} \]
and by uniform continuity the second term can be made arbitrarily small as $\varepsilon \to 0$. We can then easily show a similar inequality for every $0 < t < T$. 

Moreover, outside of $[-T, T]$, it is clear that $\overline{e}$ is still in $C^0_t L^1_x$. 

Hence we constructed a solution in $L^{\infty}_{t, loc} L^2_x \cap C^0_t L^2_{w, x}$. If we define the norm $L^2_{\mathbf{K}^{ \varepsilon}(t)}$ as in lemma \ref{lem-uniforme-convexite-global} from $x \mapsto \mathbf{K}^{\varepsilon}(t, x)$, we have that 
\[ t \mapsto \Vert u_{1/2}(t) \Vert_{L^2_{\mathbf{K}^{\varepsilon}(t)}} = \Vert \overline{e}(t) \Vert_{L^1}^{1/2} \]
is a continuous function. Let now $t \in \mathbb{R}$ be fixed, and let us show that 
\[ \Vert u_{1/2}(s) \Vert_{L^2_{\mathbf{K}^{\varepsilon}(t)}} \underset{s \to t}{\longrightarrow} \Vert u_{1/2}(t) \Vert_{L^2_{\mathbf{K}^{\varepsilon}(t)}} \]
By applying lemma \ref{lem-uniforme-convexite-global} that proves uniform convexity of this norm (for fixed $t$), knowing that $u_{1/2}$ is continuous in $L^2_{w, x}$, we will deduce that $u_{1/2}$ is even continuous in $L^2_x$ (for the strong topology of $L^2_x$). 

Yet by lemma \ref{lem-jauge-Lipschitz}, we have that 
\begin{align*} \left| j_{\mathbf{K}^{\varepsilon}(t, x)}(u_{1/2}(s, x)) - j_{\mathbf{K}^{\varepsilon}(s, x)}(u_{1/2}(s, x)) \right| &\leq d_{\mathcal{H}}(\mathbf{K}^{\varepsilon}(t, x), \mathbf{K}^{\varepsilon}(s, x)) \frac{|u_{1/2}(s, x)|}{r_{min}(\mathbf{K}^{\varepsilon}(s, x)) r_{min}(\mathbf{K}^{\varepsilon}(t, x))} \\
&\leq d_{\mathcal{H}}(\mathbf{K}^{\varepsilon}(t, x), \mathbf{K}^{\varepsilon}(s, x)) \frac{|u_{1/2}(s, x)|}{r_0^2} \end{align*}
hence 
\[ \begin{aligned}
&\int_{\mathbb{R}^n} \left| j_{\mathbf{K}^{\varepsilon}(t, x)}(u_{1/2}(s, x))^2 - j_{\mathbf{K}^{\varepsilon}(s, x)}(u_{1/2}(s, x))^2 \right|^2 ~ dx \\
&\leq \int_{\mathbb{R}^n} \frac{d_{\mathcal{H}}(\mathbf{K}^{\varepsilon}(t, x), \mathbf{K}^{\varepsilon}(s, x))}{r_0^2} |u_{1/2}(s, x)|^2 ~ dx 
\end{aligned} \]
But $u_{1/2} \in L^{\infty}_{t, loc} L^2_x$, and it is thus enough to show that $d_{\mathcal{H}}(\mathbf{K}^{\varepsilon}(t, x), \mathbf{K}^{\varepsilon}(s, x))$ converges to $0$ in $L^{\infty}_x$ as $s \to t$. But $\mathbf{K}^{0, \varepsilon}(t, x)$ is defined as the intersection of two balls whose centers and radii are $C^{\infty}_{t, x}$ functions, with uniform (but degenerating with $\varepsilon$) bounds on every derivative, and $\mathbf{K}^{\varepsilon}$ is its rescaling, hence the result. 
\end{Dem}

\begin{Dem}[of theorem \ref{theo-cheminsCdemi}]
By lemma \ref{lem-construction-intermediaire}, we construct a sequence $u_s$ for $s$ a dyadic number in $[0, 1]$, with for every $k \in \mathbb{N}$ and every $l = 0, ..., 2^k-1$ : 
\begin{multline}
\max \left( \Vert u_{(2l+1)/2^{k+1}} - u_{l/2^k} \Vert_{L^{\infty}([-T, T], L^2_x)}, ~ \Vert u_{(2l+1)/2^{k+1}} - u_{(l+1)/2^k} \Vert_{L^{\infty}([-T, T], L^2_x)} \right) \\
\leq \frac{\Vert u_{l/2^k} - u_{(l+1)/2^k} \Vert_{L^{\infty}([-T-2^{-k}, T+2^{-k}], L^2_x)}}{\sqrt{2}} + \delta_k \label{equ-construction-dyadiques} \end{multline}
for every $0 < T < k$, 
with $\delta_k = 2^{-k} \Vert u_{l/2^k} - u_{(l+1)/2^k} \Vert_{L^{\infty}([-k-2^{-k}, k+2^{-k}], L^2_x)}$. (The case where the difference vanishes is not taken into account by lemma \ref{lem-construction-intermediaire}, but then there is nothing to do to construct an intermediate solution.) 

Then, the sequence $u_s$ satisfy : 
\begin{equation} \Vert u_s - u_{s'} \Vert_{L^{\infty}([-1, 1], L^2_x)} \leq C \Vert u_0 - u_1 \Vert_{L^{\infty}([-2, 2], L^2_x)} |s - s'|^{1/2} \label{equ-Holder-dyadique} \end{equation}
for some universal constant $C > 0$, iterating the relation \ref{equ-construction-dyadiques}. 

This means that the map $s \mapsto u_s$, for now only defined on dyadic numbers and restricted to $[-1, 1]$ in time, is a $\frac{1}{2}$-Hölder functions valued in $C^0_t L^2_x$. But this space is complete, so this map admits a unique extension as a $\frac{1}{2}$-Hölder function $[0,1] \to C^0_t L^2_x$. Since the incompressible Euler equation is closed in $C^0_t L^2_x$, the whole path remains in the set of weak solutions. 

To obtain the same result for any $T$, we choose without loss of generality $T = k$ for some $k$, and then we notice that the construction for $T = 1$ on any interval $[l 2^{-k}, (l+1) 2^{-k}]$ gives Hölder regularity $\frac{1}{2}$, valued in $C^0([-k, k], L^2_x)$, then that a continuous function which is also of Hölder regularity by pieces is globally Hölder regular. 
\end{Dem}

\section{Open questions and perspectives}

We now provide some open questions related to the geometry of the set of weak solutions to the incompressible Euler equations. 

\begin{enumerate}
\item Theorem \ref{theo-cheminsCdemi} provides a Hölder regularity $\frac{1}{2}$. It is natural to ask whether this regularity is optimal. 

On the one hand, the construction presented above inherits from the procedure of De Lellis and Szekelyhidi some features that make it impossible for the path to be $C^1$, even if the convex sets are chosen differently in a more subtle manner. More precisely, the main obstruction comes from the localisation in space and time. Consider the case where $u_0 \equiv 0$ is the $0$ solution and $u_1$ is some weak solution from \cite{DeLellisSzekelyhidioriginal} with compact support in space and time. Assume by contradiction that a path following a similar convex integration procedure is of regularity $C^1$, and in particular of regularity $C^1$ at $s = 0$. The procedure ensures that, along the path, $\gamma(s)$ conserves a similar localisation property as the final data $u_1$. 

However, if we denote by $v = \partial_s \gamma(0)$, then by $C^1$ regularity at $0$ we deduce that $v$ satisfies the following equation : 
\[ \left\{ \begin{array}{l}
\partial_t v + \nabla_x q = 0 \\
\nabla_x \cdot v = 0
\end{array} \right. \]
for some $q$. This implies that $v(t) = v(0)$ is independent of $t$. But this would mean that for small $s$, $\gamma(s)$ cannot be localised in time, and this is a contradiction. 

Note that this argument does not discard the existence of $C^1$ path, but merely prevents them to be reached by a convex integration procedure that preserves localisation. 

On the other hand, if $u_0$ and $u_1$ are both classical solutions, by interpolating linearly their initial data and constructing a classical solution for every $s$ (at least locally in time, but it can be made global in dimension $n = 2$ at least) one obtains a $C^1$ path. 

Interestingly, the case of Yudovich's solutions is not clear, since the stability (ie the continuity with respect to initial data) does not follow from the proof of their unicity. 

\item We do not fix any sort of initial data in theorem \ref{theo-cheminsCdemi}, and there is no reason to hope for the initial data to be preserved in our construction if both $u_0$ and $u_1$ have the same. The main obstruction comes from the regularisation in time. Namely, mollification by convolution is difficult to consider whenever there is a boundary condition. It is open to know whether a more precise construction could prove path-connectedness of sets of weak solutions sharing the same initial data. 

\item Likewise, we avoid any boundary issue by considering the full space (although the construction is the same for instance on the torus). A construction in a setting with space boundary and boundary conditions should be of similar difficulty as point 2 above. 

\item The solutions constructed along the path do not satisfy any kind of energy inequality. Indeed, in comparison with \cite{DeLellisSzekelyhidioriginal}, we do not track the energy $|u|^2$ but a non-standard distance to a base point (which is not $0$). It seems however natural to ask whether solutions satisfying an energy inequality are also path-connected, and how to construct such paths. 

Note that this question can be combined with point 2 above, but in this case any convex integration procedure likely fails. Indeed, it was proven in \cite{BrenierDLSzWeakstrongunique} that the incompressible Euler solutions with an energy inequality enjoy (even at a measure-valued level) weak-strong uniqueness. 
\end{enumerate}

\section*{Acknowledgments}

The research of this paper has been prompted by a math reading group on convex integration organized at the École Normale Supérieure by Emmanuel Giroux, to whom the author would like to express his gratitude. The author also thanks Yann Brenier and Frédéric Rousset for helpful discussions and feedback.


\begin{thebibliography}{0}

\bibitem{BrenierDLSzWeakstrongunique} Y. Brenier, C. De~Lellis and L. Sz\'ekelyhidi Jr., Weak-strong uniqueness for measure-valued solutions, Comm. Math. Phys. {\bf 305} (2011), no.~2, 351--361. 

\bibitem{BLNWild} A.~C. Bronzi, M.~C. Lopes-Filho and H.~J. Nussenzveig~Lopes, Wild solutions for 2D incompressible ideal flow with passive tracer, Commun. Math. Sci. {\bf 13} (2015), no.~5, 1333--1343. 

\bibitem{BuckmasterOnsagerae} T. Buckmaster, Onsager's conjecture almost everywhere in time, Comm. Math. Phys. {\bf 333} (2015), no.~3, pp.~1175--1198.

\bibitem{BDLSOnsager} T. Buckmaster, C. De~Lellis and L. Sz\'ekelyhidi Jr., Dissipative Euler flows with Onsager-critical spatial regularity, Comm. Pure Appl. Math. {\bf 69} (2016), no.~9, 1613--1670. 

\bibitem{BDLISAnomalousdissip} T. Buckmaster, C. De~Lellis, P. Isett and L. Sz\'ekelyhidi Jr., Anomalous dissipation for $1/5$-H\"older Euler flows, Ann. of Math. (2) {\bf 182} (2015), no.~1, 127--172. 

\bibitem{BDLSVOnsager} T. Buckmaster, C. De~Lellis, L. Sz\'ekelyhidi Jr. and V. Vicol, Onsager's conjecture for admissible weak solutions, Comm. Pure Appl. Math. {\bf 72} (2019), no.~2, 229--274. 

\bibitem{BSV_SQG} T. Buckmaster, S. Shkoller and V.~C. Vicol, Nonuniqueness of weak solutions to the SQG equation, Comm. Pure Appl. Math. {\bf 72} (2019), no.~9, 1809--1874. 

\bibitem{BVNS} T. Buckmaster and V.~C. Vicol, Nonuniqueness of weak solutions to the Navier-Stokes equation, Ann. of Math. (2) {\bf 189} (2019), no.~1, 101--144.

\bibitem{BVreview} T. Buckmaster and V.~C. Vicol, Convex integration constructions in hydrodynamics, Bull. Amer. Math. Soc. (N.S.) {\bf 58} (2021), no.~1, 1--44.

\bibitem{CDLKisentropic} E. Chiodaroli, C. De~Lellis and O. Kreml, Global ill-posedness of the isentropic system of gas dynamics, Comm. Pure Appl. Math. {\bf 68} (2015), no.~7, 1157--1190. 

\bibitem{CMBoussinesq} E. Chiodaroli and M. Mich\'alek, Existence and non-uniqueness of global weak solutions to inviscid primitive and Boussinesq equations, Comm. Math. Phys. {\bf 353} (2017), no.~3, 1201--1216. 

\bibitem{CSzStationary} A. Choffrut and L. Sz\'ekelyhidi Jr., Weak solutions to the stationary incompressible Euler equations, SIAM J. Math. Anal. {\bf 46} (2014), no.~6, 4060--4074. 

\bibitem{CFGporous} D. C\'ordoba~Gazolaz, D. Faraco and F. Gancedo~Garc\'ia, Lack of uniqueness for weak solutions of the incompressible porous media equation, Arch. Ration. Mech. Anal. {\bf 200} (2011), no.~3, 725--74. 

\bibitem{DeLellisSzekelyhidioriginal} C. De~Lellis and L. Sz\'ekelyhidi Jr., The Euler equations as a differential inclusion, Ann. of Math. (2) {\bf 170} (2009), no.~3, 1417--1436. 

\bibitem{DeLellisSzekelyhidi2} C. De~Lellis and L. Sz\'ekelyhidi Jr., On admissibility criteria for weak solutions of the Euler equations, Arch. Ration. Mech. Anal. {\bf 195} (2010), no.~1, 225--260.  

\bibitem{DeLellisSzekelyhidi3} C. De~Lellis and L. Sz\'ekelyhidi Jr., Dissipative continuous Euler flows, Invent. Math. {\bf 193} (2013), no.~2, 377--407. 

\bibitem{DeLellisSzekelyhidi4} C. De~Lellis and L. Sz\'ekelyhidi Jr., On $h$-principle and Onsager's conjecture, Eur. Math. Soc. Newsl. No. 95 (2015), 19--24. 

\bibitem{Eulerequation1755} L. Euler, Opera Omnia, Series Secunda, 12. 

\bibitem{FLSMHD1} D. Faraco, S. Lindberg and L. Sz\'ekelyhidi Jr., Bounded solutions of ideal MHD with compact support in space-time, Arch. Ration. Mech. Anal. {\bf 239} (2021), no.~1, 51--93. 

\bibitem{FLSMHD2} D. Faraco, S. Lindberg and L. Sz\'ekelyhidi Jr., Magnetic helicity, weak solutions and relaxation of ideal MHD, Comm. Pure Appl. Math. {\bf 77} (2024), no.~4, 2387--2412. 

\bibitem{Gromovbook1986} M. Gromov, {\it Partial differential relations}, Ergeb. Math. Grenzgeb., 3. Folge, 9, Springer (1986). 

\bibitem{IsettOnsager} P. Isett, A proof of Onsager's conjecture, Ann. of Math. (2) {\bf 188} (2018), no.~3, 871--963. 

\bibitem{IsettOnsager2} P. Isett, On the endpoint regularity in Onsager's conjecture, Anal. PDE {\bf 17} (2024), no.~6, 2123--2159. 

\bibitem{KYPerona} S. Kim and B. Yan, Convex integration and infinitely many weak solutions to the Perona-Malik equation in all dimensions, SIAM J. Math. Anal. {\bf 47} (2015), no.~4, 2770--2794. 

\bibitem{Kolmogorovturbulence} A.~N. Kolmogorov, The local structure of turbulence in incompressible viscous fluid for very large Reynold's numbers, C. R. (Doklady) Acad. Sci. URSS (N.S.) {\bf 30} (1941), 301--305.

\bibitem{MajdabookEuler} A.~J. Majda and A.~L. Bertozzi, {\it Vorticity and incompressible flow}, Cambridge Texts in Applied Mathematics, 27, Cambridge Univ. Press, Cambridge, 2002. 

\bibitem{NashC1embeddings} J.~F. Nash Jr., $C^1$ isometric imbeddings, Ann. of Math. (2) {\bf 60} (1954), 383--396. 

\bibitem{Ngeostrophic} M.~D. Novack, On the weak solutions to the three-dimensional inviscid quasi-geostrophic system, SIAM J. Math. Anal. {\bf 51} (2019), no.~3, 2686--2712. 

\bibitem{Onsagerhydro} L. Onsager, Statistical hydrodynamics, Nuovo Cimento (9) {\bf 6} (1949), Supplemento, no. 2 (Convegno Internazionale di Meccanica Statistica), 279--287.

\bibitem{Schefferfirstconstr} V. Scheffer, An inviscid flow with compact support in space-time, J. Geom. Anal. {\bf 3} (1993), no.~4, 343--401. 

\bibitem{Schnirelmansecondconstr} A. Shnirelman, On the nonuniqueness of weak solution of the Euler equation, Comm. Pure Appl. Math. (50) {\bf 12} (1997), 1261--1286. 

\bibitem{ShvActivescalar} R. Shvydkoy, Convex integration for a class of active scalar equations, J. Amer. Math. Soc. {\bf 24} (2011), no.~4, 1159--1174. 

\bibitem{SzePorous} L. Sz\'ekelyhidi Jr., Relaxation of the incompressible porous media equation, Ann. Sci. \'Ec. Norm. Sup\'er. (4) {\bf 45} (2012), no.~3, 491--509. 

\bibitem{Tartarcompensatedcompactnesspde} L. Tartar, Compensated compactness and applications to partial differential equations, in {\it Nonlinear analysis and mechanics: Heriot-Watt Symposium, Vol. IV}, pp. 136--212, Res. Notes in Math., 39, Pitman, Boston, Mass.-London. 

\bibitem{Tartarcompensatedcompactnesscons} L. Tartar, The compensated compactness method applied to systems of conservation laws, in {\it Systems of nonlinear partial differential equations (Oxford, 1982)}, 263--285, NATO Adv. Sci. Inst. Ser. C: Math. Phys. Sci., 111, Reidel, Dordrecht. 

\bibitem{VillaniBourbaki} C. Villani, Paradoxe de Scheffer-Shnirelman revu sous l'angle de l'int\'egration convexe (d'apr\`es C. De Lellis et L. Sz\'ekelyhidi), Ast\'erisque No. 332 (2010), Exp. No. 1001, vii, 101--134.
\end{thebibliography}
\end{document}